%% file: 2002-2.tex
\documentclass{gtart}

\input gtoutput

\lognumber{212}
\volumenumber{6}\papernumber{2}\volumeyear{2002}
\pagenumbers{27}{58}
\accepted{25 January 2002}
\proposed{Robion Kirby}
\seconded{Ronald Stern, Ronald Fintushel}
\received{16 October 2001}
\published{29 January 2002}

\usepackage{amsmath,amssymb}

\newtheorem{thm}{Theorem}[section]
\newtheorem{prop}[thm]{Proposition}
\newtheorem{lemma}[thm]{Lemma}

\newtheorem{cor}[thm]{Corollary}
\theoremstyle{definition}
\newtheorem{defn}[thm]{Definition}
\newtheorem{rem}[thm]{Remark}

\newcommand{\arr}{{\mathbb R}}
\newcommand{\zee}{{\mathbb Z}}
\newcommand{\cee}{{\mathbb C}}
\newcommand{\cue}{{\mathbb Q}} 
\newcommand{\Fix}{{\mbox{\rm Fix}}}

\newcommand{\spinc}{{\mbox{\rm spin${}^c$}} }
\newcommand{\dbar}{\mbox{$\bar{\partial}$}}

\newcommand{\grad}{\nabla}

\newcommand{\Zset}{\zee}
\newcommand{\Tr}{{\mbox{\rm Tr}}}
\newcommand{\Sym}{{\mbox{\rm Sym}}}
\newcommand{\sgn}{{\mathrm{sgn}}}
\renewcommand{\co}{\mskip 1.5mu\colon\thinspace}

\begin{document}
\title{Torsion, TQFT, and Seiberg--Witten invariants\\of 3--manifolds}
\asciititle{Torsion, TQFT, and Seiberg-Witten invariants
of 3-manifolds}
\author{Thomas Mark}
\address{Department of Mathematics, University of California\\Berkeley,
CA  94720-3840, USA}
\email{mark@math.berkeley.edu}

\begin{abstract}
We prove a conjecture of Hutchings and Lee relating the 
Seiberg--Witten invariants of a closed 3--manifold $X$ with $b_1\geq 1$ to an 
invariant 
that ``counts'' gradient flow lines---including closed orbits---of a 
circle-valued Morse function on the manifold. The proof is based on a 
method described by Donaldson for computing the Seiberg--Witten 
invariants of 3--manifolds by making use of a ``topological quantum 
field theory,'' which makes the calculation completely explicit.
We also realize a version of the Seiberg--Witten invariant of $X$ as
the intersection number of a pair of totally real submanifolds of
a product of vortex moduli spaces on a Riemann surface
constructed from geometric data on $X$. The analogy with recent
work of Ozsv\'ath and Szab\'o suggests a generalization of a
conjecture of Salamon, who has proposed a model for the
Seiberg--Witten--Floer homology of $X$ in the case that $X$ is a
mapping torus.
\end{abstract}

\asciiabstract{We prove a conjecture of Hutchings and Lee relating the
Seiberg-Witten invariants of a closed 3-manifold X with b_1 > 0 to an
invariant that `counts' gradient flow lines--including closed
orbits--of a circle-valued Morse function on the manifold. The proof
is based on a method described by Donaldson for computing the
Seiberg-Witten invariants of 3-manifolds by making use of a
`topological quantum field theory,' which makes the calculation
completely explicit.  We also realize a version of the Seiberg-Witten
invariant of X as the intersection number of a pair of totally real
submanifolds of a product of vortex moduli spaces on a Riemann surface
constructed from geometric data on X. The analogy with recent work of
Ozsvath and Szabo suggests a generalization of a conjecture of
Salamon, who has proposed a model for the Seiberg-Witten-Floer
homology of X in the case that X is a mapping torus.}

\asciikeywords{Seiberg-Witten invariant, torsion, topological quantum field theory}

\maketitlepage

\section{Introduction}
In \cite{HL2} and \cite{HL1}, Hutchings and
Lee investigate circle-valued Morse theory for Riemannian manifolds
$X$ with first Betti number $b_1 \geq 1$.  Given a generic Morse
function $\phi\co  X\to S^1$ representing an element of infinite order in
$H^1(X;\Zset)$ and having no extrema, they determine a relationship
between the Reidemeister torsion $\tau(X,\phi)$ associated to
$\phi$, which is in general an element of the field $\cue(t)$, and the
torsion of a ``Morse complex'' $M^*$ defined over the ring
$L_\Zset$ of integer-coefficient Laurent series in a single
variable $t$. If $S$ is the inverse image of a regular value of
$\phi$ then upward gradient flow of $\phi$ induces a return map
$F\co  S\to S$ that is defined away from the descending manifolds of
the critical points of $\phi$. The two torsions $\tau(X,\phi)$
and $\tau(M^*)$ then differ by multiplication by the zeta
function $\zeta(F)$. In the case that $X$ has dimension three,
which will be our exclusive concern in this paper, the statement
reads
\begin{equation}
\tau(M^*)\zeta(F) = \tau(X,\phi),
\label{HLeqn}
\end{equation}
up to multiplication by $\pm t^k$. 
One should think of the left-hand side as ``counting'' gradient flows 
of $\phi$; $\tau(M^*)$ is concerned with gradient flows between 
critical points of $\phi$, while $\zeta(F)$, defined in terms of 
fixed points of the return map, describes the closed orbits of
$\phi$. It should be remarked that $\tau(X,\phi)\in \cue(t)$ is in
fact a polynomial if $b_1(X)>1$, and ``nearly'' so if $b_1(X) =
1$; see \cite{MT} or \cite{Turaev1} for details.

If the three--manifold $X$ is zero-surgery on a knot 
$K\subset S^3$ and $\phi$ represents a generator in $H^1(X;\Zset)$, 
the 
Reidemeister torsion $\tau(X,\phi)$
is essentially (up to a standard factor) the Alexander polynomial 
$\Delta_K$ of the knot. It has been proved by Fintushel and Stern 
\cite{FS} that 
the Seiberg--Witten invariant of $X\times S^1$, which can be 
identified with the Seiberg--Witten invariant of $X$, is also given by 
the Alexander 
polynomial (up to the same standard factor). More generally, Meng and 
Taubes \cite{MT} show that the Seiberg--Witten invariant of any closed 
three--manifold with $b_1(X)\geq 1$ can be identified with the Milnor 
torsion 
$\tau(X)$ (after summing over the action of the torsion subgroup of 
$H^2(X;\Zset)$), from which it follows that if $\mathcal S$ denotes 
the 
collection of spin${}^c$ structures on $X$,
\begin{equation}
\sum_{\alpha\in {\mathcal S}} SW(\alpha) t^{c_1(\alpha)\cdot {S} /2} 
= 
\tau(X,\phi),
\label{MTthm}
\end{equation}
up to multiplication by $\pm t^k$ (in \cite{MT} the sign is specified). 
Here $c_1(\alpha)$ denotes the first
Chern class of the complex line bundle $\det \alpha$ associated to
$\alpha$. 

These results point to the natural conjecture, made in \cite{HL1}, 
that the left-hand side of (\ref{HLeqn}) is equal to the 
Seiberg--Witten invariant of $X$---or more precisely to a combination 
of invariants as in (\ref{MTthm})---independently of the results of 
Meng and Taubes. We remark that the theorem of Meng and Taubes
announced in \cite{MT} depends on surgery formulae for
Seiberg--Witten invariants, and a complete proof of these results
has not yet appeared in the literature. The conjecture of
Hutchings and Lee gives a direct interpretation of the
Seiberg--Witten invariants in terms of geometric information,
reminiscent of Taubes's work relating Seiberg--Witten invariants
and holomorphic curves on symplectic 4--manifolds. The proof of
this conjecture is the aim of this paper; combined with the work
in \cite{HL1} and \cite{HL2} it establishes an alternate proof of
the Meng--Taubes result (for closed manifolds) that does not
depend on the surgery formulae for Seiberg--Witten invariants used
in \cite{MT} and \cite{FS}.

\begin{rem} In fact, the conjecture in \cite{HL1} is more general, 
as follows: 
Hutchings and Lee define an invariant $I\co  {\mathcal S} \to \Zset$ of 
spin${}^c$ structures based on the counting of gradient flows, which 
is conjectured to agree with the Seiberg--Witten invariant. The proof 
presented in this paper gives only an ``averaged'' version of this 
statement, ie, that the left hand side of (\ref{HLeqn}) is equal 
to the left hand side of (\ref{MTthm}). It can be seen from the 
results 
of \cite{HL1} that this averaged statement is in fact enough to 
recover the full Meng--Taubes theorem: see in particular \cite{HL1}, 
Lemma 
4.5. It may also be possible to extend the methods of this paper to 
distinguish the Seiberg--Witten invariants of spin${}^c$ structures 
whose determinant lines differ by a non-torsion element $a\in 
H^2(X;\Zset)$ with $a\cdot {S} = 0$.
\end{rem}

We also show that the ``averaged'' Seiberg--Witten invariant is equal
to the intersection number of a pair of totally real submanifolds in a
product of symmetric powers of a slice for $\phi$.  This is a
situation strongly analogous to that considered by Ozsv\'ath and
Szab\'o in \cite{OS1} and \cite{OS2}, and one might hope
to define a Floer-type homology theory along the lines of that work. 
Such a construction would suggest a generalization of a conjecture of
Salamon, namely that the Seiberg--Witten--Floer homology of $X$ agrees
with this new homology (which is a ``classical'' Floer homology in the
case that $X$ is a mapping torus---see \cite{S}).

\section{Statement of results}

Before stating our main theorems, we need to recall a few definitions
and introduce some notation.  First is the notion of the torsion of 
an acyclic
chain complex; basic references for this material include \cite{Milnor}
and \cite{Turaev1}.

\subsection{Torsion}
By a {\it volume}
 $\omega$ for a vector space $W$ of dimension $n$ we mean a
choice of nonzero element $\omega\in\Lambda^n W$. Let $0\to V'\to
V\to V''\to 0$ be an exact sequence of finite-dimensional vector
spaces over a field $k$. For volumes $\omega'$ on $V'$ and
$\omega''$ on $V''$, the induced volume on $V$ will be written
$\omega'\omega''$; if $\omega_1$, $\omega_2$ are two volume
elements for $V$, then we can write $\omega_1 = c\omega_2$ for
some nonzero element $c\in k$ and by way of shorthand, write $c =
\omega_1/\omega_2$. More generally, let $\{C_i\}_{i=0}^n$ be a
complex of vector spaces with differential $\partial\co  C_i\to
C_{i-1}$, and let us assume that $C_*$ is acyclic, ie,
$H_*(C_*)=0$. Suppose that each $C_i$ comes equipped with a
volume element $\omega_i$, and choose volumes $\nu_i$ arbitrarily
on each image $\partial C_i$, $i=2,\ldots,n-1$. From the exact sequence
\[
0\to C_n \to C_{n-1} \to \partial C_{n-1}\to 0
\]
define $\tau_{n-1} = \omega_n\nu_{n-1}/\omega_{n-1}$. For $i= 
2,\ldots, n-2$ use 
the exact sequence
\[
0\to \partial C_{i+1}\to C_i\to \partial C_i\to 0
\]
to define $\tau_i = \nu_{i+1}\nu_i/\omega_i$. Finally, from
\[
0\to \partial C_2\to C_1\to C_0\to 0
\]
define $\tau_1 = \nu_2\omega_0/\omega_1$. We then define
the {\it torsion} $\tau(C_*, \{\omega_i\}) \in k\setminus\{0\}$
of the (volumed) complex $C_*$ to be:
\begin{equation}
\tau(C_*) = \prod_{i=1}^{n-1} \tau_i^{(-1)^{i+1}}
\label{torsiondef}
\end{equation}
It can be seen that this definition does not depend on the choice of 
$\nu_i$. 
Note that in the case that our complex consists of just two vector 
spaces, 
\[
C_* = 0\to C_i\stackrel{\partial}{\longrightarrow} C_{i-1}\to 0,
\]
we have that $\tau(C)= \det(\partial)^{(-1)^i}$. We extend the
definition of $\tau(C_*)$ to non-acyclic complexes by setting
$\tau(C_*) = 0$ in this case.

As a slight generalization, we can allow the chain groups $C_i$ to be 
finitely generated free 
modules over an integral domain $K$ with fixed ordered bases rather 
than vector 
spaces with fixed volume elements, as 
follows. Write $Q(K)$ for the field of fractions of $K$, then form 
the complex of vector spaces $Q(K)\otimes_K C_i$. The bases for 
the $C_i$ naturally give rise to bases, and hence volumes, for
$Q(K)\otimes_K C_i$. We understand the torsion of the complex of
$K$--modules $C_i$ to be the torsion of this latter complex, and
it is therefore a nonzero element of the field $Q(K)$.

Let $X$ be a connected, compact, oriented smooth manifold with a
given CW decomposition. Following \cite{Turaev1}, suppose
$\varphi\co  \zee[H_1(X;\zee)]\to K$ is a ring homomorphism into an
integral domain $K$. The universial abelian cover $\tilde{X}$ has
a natural CW decomposition lifting the given one on $X$, and the
action of the deck transformation group $H_1(X;\zee)$ naturally
gives the cell chain complex $C_*(\tilde{X})$ the structure of a
$\zee[H_1(X;\zee)]$--module. As such, $C_i(\tilde{X})$ is free of
rank equal to the number of $i$--cells of $X$. We can then form
the twisted complex $C_*^\varphi(\tilde{X}) = K\otimes_\varphi
C_*(\tilde{X})$ of $K$--modules. We choose a sequence $e$ of cells
of $\tilde{X}$ such that over each cell of $X$ there is exactly
one element of $e$, called a {\it base} {\it sequence}; this gives
a basis of $C_*^\varphi(\tilde{X})$ over $K$ and allows us to
form the torsion $\tau_\varphi(X,e)\in Q(K)$ relative to this
basis. Note that the torsion $\tau_\varphi(X,e')$ arising from a
different choice $e'$ of base sequence stands in the relationship
$\tau_\varphi(X,e) = \pm\varphi(h)\tau_\varphi(X,e')$ for some
$h\in H_1(X;\zee)$ (here, as is standard practice, we write the
group operation in $H_1(X;\zee)$ multiplicatively when dealing
with elements of $\zee[H_1(X;\zee)]$). The set of all torsions
arising from all such choices of $e$ is ``the'' torsion of $X$
associated to $\varphi$ and is denoted $\tau_\varphi(X)$.

We are now in a position to define the torsions we will need.

\begin{defn}(1)\qua For $X$ a smooth manifold as above with $b_1(X)\geq 
1$, 
let $\phi\co  X\to S^1$ be a map 
representing an element $[\phi]$ of infinite order in $H^1(X;\zee)$. 
Let $C$ 
be the infinite cyclic group generated by the formal variable $t$, 
and let
$\varphi_1\co  \zee[H_1(X;\zee)]\to \zee[C]$ be the map induced 
by the 
homomorphism $H_1(X;\zee)\to C$, $a \mapsto t^{\langle 
[\phi],a\rangle}$. 
Then the {\em Reidemeister torsion} $\tau(X,\phi)$ of $X$ associated
to $\phi$ is defined to be the torsion $\tau_{\varphi_1}(X)$.

(2)\qua Write $H$ for the quotient of $H_1(X;\zee)$ by its torsion 
subgroup, and let $\varphi_2\co  \zee[H_1(X;\zee)]\to\zee[H]$ be the map 
induced by the projection $H_1(X;\zee)\to H$. The {\em Milnor 
torsion} $\tau(X)$ is defined to be $\tau_{\varphi_2}(X)$.
\label{torsiondefn} 
\end{defn}

\begin{rem}(1)\qua Some authors use the term {\it Reidemeister torsion} 
to 
refer to the torsion $\tau_\varphi(X)$ for arbitrary $\varphi$; and other terms, 
eg, 
Reidemeister--Franz--DeRham torsion, are also in use.

(2)\qua The torsions in Definition \ref{torsiondefn} are defined for 
manifolds $X$ of arbitrary dimension, with or without boundary. We 
will be concerned only with the case that $X$ is a closed manifold of 
dimension 3 with $b_1(X)\geq 1$.  In the case $b_1(X)>1$, work of 
Turaev
\cite{Turaev1} shows that $\tau(X)$ and $\tau(X,\phi)$, naturally
subsets of $\cue(H)$ and $\cue(t)$, are actually subsets of $\zee[H]$
and $\zee[t, t^{-1}]$. Furthermore, if $b_1(X)=1$ and $[\phi]\in
H^1(X;\zee)$ is a generator, then $\tau(X) = \tau(X,\phi)$ 
and $(t-1)^2\tau(X)\in\zee[t,t^{-1}]$. Rather than thinking of
torsion as a set of elements in a field we normally identify it
with a representative ``defined up to multiplication by $\pm
t^k$'' or similar, since by the description above any two
representatives of the torsion differ by some element of the
group ($C$ or $H$) under consideration.
\end{rem}

\subsection{$S^1$--Valued Morse Theory}
\label{morsesec}

We review the results of Hutchings and Lee that motivate our theorems. 
As in the introduction, let $X$ be a smooth closed oriented 3--manifold
having $b_1(X)\geq 1$ and let $\phi\co X\to S^1$ be a smooth Morse
function.  We assume (1) $\phi$ represents an indivisible element of 
infinite
order in $H^1(X,\zee)$; (2) $\phi$ has no critical points of index
0 or 3; and (3) the gradient flow of $\phi$ with respect to a
Riemannian metric on $X$ is Morse--Smale. Such functions always exist
given our assumptions on $X$.

Given such a Morse function $\phi$, fix a smooth level set $S$ for
$\phi$. Upward gradient flow defines a return map $F\co  S\to S$
away from the descending manifolds of the critical points of
$\phi$. The {\it zeta function} of $F$ is defined by the series
\[
\zeta(F) = \exp \left(\sum_{k\geq 1} \Fix(F^k)\frac{t^k}{k}\right)
\]
where $\Fix(F^k)$ denotes the number of fixed points (counted with
sign in the usual way) of the $k$-th iterate of $F$.  One should think
of $\zeta(F)$ as keeping track of the number of closed orbits of
$\phi$ as well as the ``degree'' of those orbits.  For future
reference we note that if $h\co S\to S$ is a diffeomorphism of a surface 
$S$ then 
\begin{equation}
\zeta(h) = \sum_k L(h^{(k)})t^k
\label{zetasym}
\end{equation}
where $L(h^{(k)})$ is the Lefschetz number of the induced map on the
$k$-th symmetric power of $S$ (see \cite{S}, \cite{IP}).

We now introduce a Morse complex that can be used to keep track of
gradient flow lines between critical points of $\phi$. Write 
$L_{\zee}$ for the ring of Laurent series in 
the variable $t$, and let $M^i$ denote the free $L_{\zee}$--module 
generated by 
the index-$i$ critical points of $\phi$.  The differential $d_M\co M^i\to 
M^{i+1}$ is defined to be
\begin{equation}
d_Mx_\mu = \sum_\nu a_{\mu\nu}(t) y_\nu
\label{diffdef}
\end{equation}
where $x_\mu$ is an index-$i$ critical point, $\{y_\nu\}$
is the set of index-$(i+1)$ critical points, and $a_{\mu\nu}(t)$ is a
series in $t$ whose coefficient of $t^n$ is defined to be the number
of gradient flow lines of $\phi$ connecting $x_\mu$ with $y_\nu$ that
cross $S$ $n$ times.  Here we count the gradient flows with sign
determined by orientations on the ascending and descending manifolds
of the critical points; see \cite{HL1} for more details.

\begin{thm}[Hutchings--Lee] In this situation, the relation 
(\ref{HLeqn}) holds up to multiplication by $\pm t^k$.
\label{HLthm}
\end{thm}

\subsection{Results}

The main result of this work is that the left hand side of
(\ref{HLeqn}) is equal to the left hand side of (\ref{MTthm}),
without using the results of \cite{MT}. Hence the current work,
together with that of Hutchings and Lee, gives an alternative
proof of the theorem of Meng and Taubes in \cite{MT}.

Our proof of this fact is based on ideas of 
Donaldson for computing the Seiberg--Witten invariants of 3--manifolds. 
We outline Donaldson's construction here; see Section \ref{tqftsec} 
below for 
more details. Given $\phi\co X\to S^1$ a generic Morse 
function as above and $S$ the inverse image of a regular value, 
let $W = X\setminus nbd(S)$ be the complement of a small 
neighborhood of $S$. Then $W$ is a cobordism between two copies 
of $S$ (since we assumed $\phi$ has no extrema---note we may 
also assume $S$ is connected). Note that two spin${}^c$ structures on 
$X$ that 
differ by an element $a\in H^2(X; \zee)$ with $a([S]) 
= 0$ 
restrict to the same spin${}^c$ structure on $W$, in particular, 
spin${}^c$ 
structures $\sigma$ on $W$ are determined by their degree $m = \langle
c_1(\sigma), S\rangle$. Note that the degree of a \spinc structure is 
always even.

Now, a solution of the Seiberg--Witten equations on $W$ restricts to a 
solution of the {\it vortex equations} on $S$ at each end of 
$W$ (more accurately, we should complete $W$ by adding infinite 
tubes $S\times (-\infty, 0]$, $S\times [0,\infty)$ to each 
end, and consider the limit of a finite-energy solution on this completed 
space)---see \cite{D2}, \cite{MOY} for example. These equations have 
been extensively studied, and it is known that the moduli space of 
solutions to the vortex equations on $S$ can be identified with 
a symmetric power $\mathrm{Sym}^n S$ of $S$ itself: see \cite{B}, \cite{JT}. Donaldson
uses the restriction maps on the Seiberg--Witten moduli space of $W$ to
obtain a self-map $\kappa_n$ of the cohomology of
$\mbox{Sym}^{n}S$, where $n$ is defined by $n=g(S)-1-\frac{1}{2}|m|$ if 
$b_1(X)>1$ and $n= g(S)-1 + \frac{1}{2}m$ if $b_1(X) = 1$ (here
$g(S)$ is the genus of the orientable surface $S$). The alternating trace
$\Tr\,\kappa_n$ is identified as the sum of Seiberg--Witten
invariants of spin${}^c$ structures on $X$ that restrict to the
given spin${}^c$ structure on $W$---that is, the coefficient of
$t^n$ on the left hand side of (\ref{MTthm}). For a precise statement,
see Theorem \ref{tracethm}.

Our main result is the following. 

\begin{thm} Let $X$ be a Riemannian 3--manifold with $b_1(X)\geq 1$, 
and fix an integer $n\geq 0$ as above.  Then we have
\begin{equation}
\Tr\,\kappa_n = [\tau(M^*)\,\zeta(F)]_{n},
\end{equation}
where $\tau(M^*)$ is represented by $t^N\det(d_M)$, and $N$ is the 
number of index 1 critical points of $\phi$. Here $\Tr$ denotes the 
alternating trace and $[\,\cdot\,]_n$ denotes the coefficient of $t^n$ 
of the polynomial enclosed in brackets.
\label{mainthm}
\end{thm}

This fact immediately implies the conjecture of Hutchings and Lee. 
Furthermore, we will make the following observation:

\begin{thm}
There is a smooth connected representative $S$ for the Poincar\'e dual
of $[\phi]\in H^1(X;\zee)$ such that $\Tr\, \kappa_n$ is
given by the intersection number of a pair of totally real embedded
submanifolds in $\mathrm{Sym}^{n+N}S \times\mathrm{Sym}^{n+N}S$.
\label{intthm}
\end{thm}

This may be the first step in defining a Lagrangian-type Floer 
homology theory parallel to that of Ozsv\'ath and Szab\'o, one whose 
Euler 
characteristic is {\it a priori} a combination of Seiberg--Witten 
invariants. In the case that $X$ is a mapping torus, a program along 
these lines has been initiated by Salamon \cite{S}. In this 
case the two totally real submanifolds in Theorem \ref{intthm} reduce 
to the diagonal and the graph of a symplectomorphism of
$\mathrm{Sym}^n S$ determined by the monodromy of the mapping
torus, both of which are in fact Lagrangian.

The remainder of the paper is organized as follows: Section 3 gives a 
brief overview of some elements of Seiberg--Witten theory and the 
dimensional reduction we will make use of, and Section 4 gives a few 
more details on this reduction and describes the TQFT we use to 
compute Seiberg--Witten invariants. Section 5 proves a theorem that 
gives a means of calculating as though a general cobordism coming from 
an $S^1$--valued Morse function of the kind we are considering posessed a 
naturally-defined monodromy map; Section 6 collects a few other 
technical results of a calculational nature, the proof of one of 
which is the content of Section 9. In Section 7 we prove Theorem 
\ref{mainthm} by a calculation that is fairly involved but is not 
essentially difficult, thanks to the tools provided by the TQFT. 
Section 8 proves Theorem \ref{intthm}.

\section{Review of Seiberg--Witten theory}

We begin with an outline of some aspects of Seiberg--Witten 
theory for a 3--manifolds. Recall that a \spinc structure on 
a 3--manifold $X$ is a lift of the oriented orthogonal frame bundle 
of $X$ to a principal $\spinc(3)$--bundle $\sigma$. There are two 
 representations of $\spinc(3) = \mbox{Spin}(3)\times 
U(1)/\pm1 = SU(2)\times U(1)/\pm 1$ that will interest us, namely
the spin representation $\spinc(3)\to SU(2)$ and also the
projection $\spinc(3)\to U(1)$ given by $[g, e^{i\theta}]\mapsto
e^{2i\theta}$. For a \spinc structure $\sigma$ the first of these
gives rise to the associated {\it spinor bundle} $W$ which is a
hermitian 2--plane bundle, and the second to the {\it determinant
line bundle} $L\cong \wedge^2 W$. We define $c_1(\sigma) :=
c_1(L)$. The Levi--Civita connection on $X$ together with a choice
of hermitian connection $A$ on $L^{1/2}$ gives rise to a
hermitian connection on $W$ that is compatible with the action of
Clifford multiplication $c\co  T^*_\cee X\to \mbox{End}_0 W$=
\{traceless endomorphisms of $W$\}, and thence to a Dirac
operator $D_A\co  \Gamma(W)\to \Gamma(W)$.

The {\it Seiberg--Witten equations} are equations for a pair 
$(A,\psi)\in{\cal A}(L)\times \Gamma(W)$ where ${\cal A}(L)$ denotes 
the space of hermitian connections on $L^{1/2}$, and read:
\begin{equation}
\begin{array}{rcl} D_A \psi &\,\,=\,\,& 0 \\
c(\star F_A + i\star\mu) &=& \psi\otimes\psi^* - \frac{1}{2}|\psi|^2
\end{array}
\label{sweqns}
\end{equation}
Here $\mu\in\Omega^2(X)$ is a closed form used as a perturbation; if
$b_1(X)>1$ we may choose $\mu$ as small as we like.

On a closed oriented 3--manifold the {\it Seiberg--Witten moduli space} is 
the set of $L^{2,2}$ solutions to the 
above equations modulo the action of the gauge group ${\cal G}= 
L^{2,3}(X;S^1)$, which acts on connections by conjugation and on 
spinors by complex multiplication. For generic choice of 
perturbation $\mu$ the moduli space ${\cal M}_\sigma$ is a compact 
zero--dimensional manifold that is smoothly cut out by its defining 
equations (if $b_1(X)>0$). There is a way to orient ${\cal M}_\sigma$ 
using a so-called homology orientation of $X$, and the {\it 
Seiberg--Witten invariant} of $X$ in the \spinc structure $\sigma$ is 
defined to be the signed count of points of ${\cal M}_\sigma$. One 
can show that if $b_1(X)>1$ then the resulting number is independent 
of all choices involved and depends only on $X$ (with its orientation); 
while if $b_1(X) = 
1$ there is a slight complication: in this case we need to make a
choice of generator $o$ for the free part of $H^1(X;\zee)$ and require
that $\langle[\mu]\cup o, [X]\rangle > \pi \langle c_1(\sigma)\cup o,
[X]\rangle$.

Suppose now that rather than a closed manifold, $X$ is isometric to a 
product $\Sigma\times \arr$ for some Riemann surface $\Sigma$. If $t$ 
is the coordinate in the $\arr$ direction, then Clifford 
multiplication by $dt$ is an automorphism of square $-1$ of $W$ and 
therefore splits $W$ into eigen-bundles $E$ and $F$ on which $dt$ 
acts 
as multiplication by $-i$ and $i$, respectively. In fact 
$F = K^{-1}E$ where $K$ is the canonical bundle of $\Sigma$, and 
$2E-K = L$, the determinant line of $\sigma$. Writing a section 
$\psi$ 
of $W$ as $(\alpha,\beta)\in\Gamma(E\oplus K^{-1}E)$, we can express 
the Dirac operator in this decomposition as:
\[
D_A\psi = \left(\begin{array}{cc} -i\frac{\partial}{\partial t} & 
\bar{\partial}_{B,J}^* \\ \dbar_{B,J} & i\frac{\partial}{\partial t} 
\end{array}\right)\left( \begin{array}{c} \alpha \\ \beta 
\end{array}\right)
\]
Here we have fixed a spin structure (with connection) $K^{1/2}$ on 
$\Sigma$ 
and noted 
that the choice of a connection $A$ on $L^{1/2} = E-K^{1/2}$ is 
equivalent 
to a choice of connection $B$ on $E$. The metric on 
$\Sigma\times\arr$ induces a complex structure $J$ and area form 
$\omega_\Sigma$ on $\Sigma$. Then $\dbar_{B,J}$ is the associated $\dbar$
operator on sections of $E$ with adjoint operator
$\dbar_{B,J}^*$.

The 2--forms $\Omega^2_\cee(\Sigma\times\arr)$ split as 
$\Omega^{1,1}(\Sigma)
\oplus 
[(\Omega^{1,0}(\Sigma)\oplus\Omega^{0,1}(\Sigma))\otimes\Omega^1_\cee(\arr)]$,
and we will write a form $\nu$ as $\Lambda\nu\cdot\omega_\Sigma +
\nu^{1,0} dt + \nu^{0,1} dt$ in this splitting. Thus $\Lambda
\nu$ is a complex function on $\Sigma\times\arr$, while
$\nu^{1,0}$ and $\nu^{0,1}$ are 1--forms on $\Sigma$. With these
conventions, the Seiberg--Witten equations become
\begin{equation}
\begin{array}{rcl} i\dot{\alpha} &\,\,=\,\,& \dbar^*_{B,J}\beta \\
i\dot{\beta} &=& -\dbar_{B,J} \alpha \\
 2\Lambda F_B- \Lambda F_K + 2i\Lambda\mu &=& 
i(|\alpha|^2-|\beta|^2) \\
(2F_B - F_K)^{1,0} + 2i\mu^{1,0} &=& \alpha\otimes 
\bar{\beta} 
\end{array}
\label{redsweqns}
\end{equation}
One can show 
that for a finite-energy solution either $\alpha$ or $\beta$ must 
identically vanish; apparently this implies any such solution is 
constant, and the above system of equations descends to $\Sigma$ when 
written in temporal gauge (ie, so the connection has no $dt$ component). 
The above equations (with
$\beta=0$) therefore reduce to the {\it vortex equations} in $E$,
which are for a pair $(B, \alpha)\in {\cal A}(E)\times \Gamma(E)$ and
read
\begin{eqnarray}
\dbar_{B,J} \alpha &=& 0 \label{vortex1}\\
i\star F_B + \frac{1}{2}|\alpha|^2 &=& \tau \label{vortex2}
\end{eqnarray}
where $\tau$ is a function on $\Sigma$ satisfying $\int \tau > 
2\pi\deg(E)$ and 
incorporates the curvature $F_K$ and perturbation above. These 
equations are well-understood, and it is known that the space of 
solutions to the vortex equations modulo $\mbox{Map}(\Sigma, S^1)$ 
is isomorphic to the space of solutions $(B,\alpha)$ of the single 
equation
\[
\dbar_{B,J} \alpha = 0
\]
modulo the action of $\mbox{Map}(\Sigma, \cee^*)$. The latter is 
naturally identified with the space of divisors of degree $d= 
\deg(E)$ on $\Sigma$ via the zeros of $\alpha$, and forms a K\"ahler 
manifold isomorphic to the $d$-th symmetric power $\mathrm{Sym}^d
\Sigma$, which for brevity we will abbreviate as $\Sigma^{(d)}$
from now on. We write ${\cal M}_d(\Sigma,J)$ (or simply ${\cal
M}(\Sigma)$) for the moduli space of vortices in a bundle $E$ of
degree $d$ on $\Sigma$.

The situation for $\alpha \equiv 0$ is analogous to 
the above: in this case $\beta$ satisfies $\dbar^*_{B,J}\beta = 0$ so 
that $\star_2\beta$ is a holomorphic section of $K\otimes E^*$. 
Replacing $\beta$ by $\star_2\beta$ shows that the Seiberg--Witten
equations reduce to the vortex equations in the bundle $K\otimes E^*$,
giving a moduli space isomorphic to $\Sigma^{(2g-2-d)}$.

\section{A TQFT for Seiberg--Witten invariants}
\label{tqftsec}

In this section we describe Donaldson's ``topological quantum field 
theory'' for computing the Seiberg--Witten invariants. Suppose $W$ is 
a 
cobordism between two Riemann surfaces $S_-$ and $S_+$. We complete 
$W$ by adding tubes $S_\pm\times [0,\infty)$ to the boundaries and 
endow the completed manifold $\hat{W}$ with a Riemannian metric that 
is a product on the ends. By considering finite-energy solutions to 
the Seiberg--Witten equations on $\hat{W}$ in some \spinc structure
$\sigma$, we can produce a Fredholm problem and show that such
solutions must approach solutions to the vortex equations on $S_\pm$. 
Following a solution to its limiting values, we obtain smooth maps
between moduli spaces, $\rho_\pm\co  {\cal M}(\hat{W})\to {\cal
M}(S_\pm)$.  Thus we can form
\begin{eqnarray*}
\kappa_\sigma = (\rho_-\otimes\rho_+)_*[{\cal M}(\hat{W})]&\in&
H_*({\cal M}(S_-)) \otimes H_*({\cal M}(S_+)) \\
&\cong& \hom (H^*({\cal M}(S_-)),H^*({\cal M}(S_+))).
\end{eqnarray*}
Here we use Poincar\'e duality and work with rational coefficients.

This is the basis for our ``TQFT:'' to a surface $S$ we associate the
cohomology of the moduli space ${\cal M}(S)$, and to a
cobordism $W$ between $S_-$ and $S_+$ we assign the homomorphism 
$\kappa_\sigma$:
\begin{eqnarray*}
S&\longmapsto & V_S = H^*({\cal M}(S))\\
W&\longmapsto & \kappa_\sigma\co  V_{S_-}\to V_{S_+}
\end{eqnarray*}
In the sequel we will be interested only in cobordisms $W$ that
satisify the topological assumption $H_1(W,\partial W) = \zee$.  Under
this asssumption, gluing theory for Seiberg--Witten solutions provides a
proof of the central property of TQFTs, namely that if $W_1$ and $W_2$
are composable cobordisms then $\kappa_{W_1\cup W_2} =
\kappa_{W_2}\circ\kappa_{W_1}$.

If $X$ is a closed oriented 3--manifold with $b_1(X)>0$ then the above
constructions can be used to calculate the Seiberg--Witten invariants
of $X$, as seen in \cite{D1}.  We now describe the procedure involved. 
Begin with a Morse function $\phi\co  X\to
S^1$ as in the introduction, and cut $X$ along the level set $S$
to produce a cobordism $W$ between two copies of $S$, which come
with an identification or ``gluing map'' $\partial_- W\to \partial_+
W$. Write $g$ for the genus of $S$. The cases $b_1(X)>1$ and $b_1(X)=1$ are
slightly different and we consider them separately.

Suppose $b_1(X)>1$, so the perturbation $\mu$ in (\ref{sweqns}) can
be taken to be small.  Consider the constant solutions to the equations
(\ref{redsweqns}) on the ends of $\hat{W}$, or equivalently the
possible values of $\rho_\pm$.  If $\beta \equiv 0$ then $\alpha$ is a
holomorphic section of $E$ and so the existence of a nonvanishing
solution requires $\deg(E)\geq 0$. Since $\mu$ is small,
integrating the third equation in (\ref{redsweqns}) tells us that
$2E-K$ is nonpositive. Hence existence of nonvanishing solutions
requires $0\leq \deg(E)\leq \frac{1}{2}\deg(K) = g-1$. If
$\alpha\equiv 0$, then $\star_2\beta$ is a holomorphic section of
$K-E$ so to avoid triviality we must have $0\leq\deg(K)-\deg(E)$,
ie, $\deg(E)\leq 2g-2$. On the other hand, integrating the
third Seiberg--Witten equation tells us that $2E-K$ is
nonnegative, so that $\deg(E)\geq g-1$. To summarize we have
shown that constant solutions to the Seiberg--Witten equations on
the ends of $\hat{W}$ in a \spinc structure $\sigma$ are just the
vortices on $S$ (with the finite-energy hypothesis). If
$\det(\sigma) = L$ a necessary condition for the existence of
such solutions is $-2g+2\leq \deg(L) \leq 2g-2$ (recall $L =
2E-K$ so in particular $L$ is even). If this condition is
satisfied than the moduli space on each end is isomorphic to
${\cal M}_n(S) \cong S^{(n)}$ where $n =
g-1-\frac{1}{2}|\deg(L)|$. Note that by suitable choice of
perturbation $\mu$ we can eliminate the ``reducible'' solutions,
ie, those with $\alpha \equiv 0 \equiv \beta$, which otherwise
may occur at the extremes of our range of values for $\deg(L)$.

Now assume $b_1(X) = 1$.  Integrating the third equation in
(\ref{redsweqns}) shows
\[
\langle c_1(\sigma), S\rangle - \frac{1}{\pi}\langle
[\mu],S\rangle = \frac{1}{2\pi}\int_S |\beta|^2- |\alpha|^2.
\]
The left hand side of this is negative by our assumption on $\mu$, and
we know that either $\alpha\equiv 0$ or $\beta \equiv 0$.  The first
of these possibilities gives a contradiction; hence $\beta\equiv 0$
and the system (\ref{redsweqns}) reduces to the vortex equations in $E$
over $S$. Existence of nontrivial solutions therefore requires
$\deg(E)\geq 0$, ie, $\deg(L)\geq 2-2g(S)$. Thus the moduli
space on each end of $\hat{W}$ is isomorphic to ${\cal
M}_n(S) \cong S^{(n)}$, where $n = \deg(E) =
g-1+\frac{1}{2}\deg(L)$ and $\deg(L)$ is any even integer at
least $2-2g(S)$.

\begin{thm}[Donaldson] Let $X$, $\sigma$, $\phi$, $S$, and $W$ be
as above. Write $\langle c_1(\sigma),[S]\rangle = m$ and define
either $n = g(S)-1 -\frac{1}{2}|m|$ or $n = g(S) - 1
+\frac{1}{2}m$ depending whether $b_1(X)>1$ or $b_1(X) = 1$. Then
if $n\geq 0$,
\begin{equation}
\Tr \,\kappa_\sigma = \sum_{\tilde{\sigma}\in{\cal S}_m}
SW(\tilde{\sigma})
\label{traceeqn}
\end{equation}
where ${\cal S}_m$ denotes the set of \spinc structures
$\tilde{\sigma}$ on $X$ such that $\langle
c_1(\tilde{\sigma}),[S]\rangle = m$. If $n<0$ then the right hand
side of (\ref{traceeqn}) vanishes. Here $\Tr$ denotes the graded
trace.
\label{tracethm}
\end{thm}

Note that with $n$ as in the theorem, $\kappa_\sigma$ is a linear map 
\[
\kappa_\sigma\co  H^*(S^{(n)})\to H^*(S^{(n)});
\]
as the trace of $\kappa_\sigma$ computes a sum of Seiberg--Witten 
invariants rather than just $SW(\sigma)$, we use the notation 
$\kappa_n$ rather than $\kappa_\sigma$.

Since $\kappa_n$ obeys the composition law, in order to determine
the map corresponding to $W$ we need only determine the map generated
by elementary cobordisms, ie, those consisting of a single 1-- or
2--handle addition (we need not consider 0-- or 3--handles by our
assumption on $\phi$).  In \cite{D1}, Donaldson uses an elegant
algebraic argument to determine these elementary homomorphisms.  To
state the result, recall that the cohomology of the $n$-th symmetric
power $S^{(n)}$ of a Riemann surface $S$ is given over $\zee$,
$\cue$, $\arr$, or $\cee$ by
\begin{equation}
H^*(S^{(n)}) = \bigoplus_{i=0}^n \Lambda^iH^1(S)\otimes
\Sym^{n-i}(H^0(S)\oplus H^2(S)).
\label{symprodcohom}
\end{equation}
Suppose that $W$ is an elementary cobordism connecting two surfaces
$\Sigma_g$ and $\Sigma_{g+1}$. Thus there is a unique critical point (of index 1) of the 
height function 
$h\co  W\to \arr$, and the ascending manifold of this critical point 
intersects $\Sigma_{g+1}$ in an essential simple closed curve that we will
denote by $c$.

Now, $c$ obviously bounds a disk $D\subset W$; the 
Poincar\'e--Lefschetz dual 
of $[D]\in H_2(W,\partial W)$ is a 1--cocycle that we will denote 
$\eta_0\in H^1(W)$. It is easy to check that $\eta_0$ is in the kernel 
of the restriction $r_1\co H^1(W)\to H^1(\Sigma_g)$, so we may complete 
$\eta_0$ to a basis $\eta_0,\eta_1,\ldots,\eta_{2g}$ of $H^1(W)$ with the 
property that
$\xi_1:=r_1(\eta_1),\ldots,\hspace{1ex}\xi_{2g}:=r_1(\eta_{2g})$
form a basis for $H^1(\Sigma_g)$. Since the restriction
$r_2\co H^1(W)\to H^1(\Sigma_{g+1})$ is injective, we know
$\bar{\xi}_0:=r_2(\eta_0), \ldots,\hspace{1ex} \bar{\xi}_{2g}:=
r_2(\eta_{2g})$ are linearly independent; note that $r_2(\eta_0)$
is just $c^*$, the Poincar\'e dual of $c$.

The choice of basis $\eta_j$ with its restrictions $\xi_j$, 
$\bar{\xi}_j$ gives rise to an inclusion $i\co H^1(\Sigma_g)\to 
H^1(\Sigma_{g+1})$ in the obvious way, namely $i(\xi_j) = 
\bar{\xi_j}$.  One may check that this map is independent of the 
choice of basis $\{\eta_j\}$ for $H^1(W)$ having $\eta_0$ as above.  
From the decomposition (\ref{symprodcohom}), we can extend $i$ to an 
inclusion $i\co  H^*(\Sigma_g^{(n)})\hookrightarrow 
H^*(\Sigma_{g+1}^{(n)})$.  Having produced this inclusion, we now 
proceed to suppress it from the notation, in particular in the 
following theorem.

\begin{thm}[Donaldson] In this situation, and with $\sigma$ and $n$ as
previously, the map $\kappa_n$ corresponding to the elementary cobordism
$W$ is given by
\[
\kappa_n(\alpha) = c^*\wedge\alpha.
\]
If $\bar{W}$ is the ``opposite'' cobordism between $\Sigma_{g+1}$ and
$\Sigma_g$, the corresponding $\kappa_n$ is given by the contraction
\[
\kappa_n(\beta) = \iota_{c^*}\beta,
\]
\label{wedgethm}
where contraction is defined using the intersection pairing on
$H^1(\Sigma_{g+1})$.
\end{thm}

This result makes the calculation of Seiberg--Witten invariants
completely explicit, as we see in the next few sections.

\section{Standardization of $X$}
\label{stdsec}

We now return to the situation of the introduction: namely, we 
consider a closed 3--manifold $X$ having $b_1(X)\geq 1$, with its 
circle-valued Morse function $\phi\co  X\to S^1$ having no critical 
points of index $0$ or $3$, and $N$ critical points of each index $1$ 
and $2$. We want to 
show how to identify $X$ with a ``standard'' manifold $M(g, N, h)$ 
that depends only on $N$ and a diffeomorphism $h$ of a Riemann 
surface of genus $g+N$. This standard manifold will be obtained from two
``compression bodies,'' ie, cobordisms between surfaces
incorporating handle additions of all the same index. Two copies
of the same compression body can be glued together along their
smaller-genus boundary by the identity map, then by a
``monodromy'' diffeomorphism of the other boundary component to
produce a more interesting 3--manifold. Such a manifold lends
itself well to analysis using the TQFT from the previous section,
as the interaction between the curves $c$ corresponding to each
handle is completely controlled by the monodromy. We now will
show that every closed oriented 3--manifold $X$ having $b_1(X)>0$
can be realized as such a glued-up union of compression bodies.

To begin with, we fix a closed oriented genus 0 surface
$\Sigma_0$ (that is, a standard 2--sphere) with an
orientation-preserving embedding $\psi_{0,0}\co S^0\times D^2\to
\Sigma_0$. Here we write $D^n = \{x\in\arr^n||x|<1\}$ for the
unit disk in $\arr^n$. There is a standard way to perform surgery
on the image of $\psi_{0,0}$ (see \cite{milnor2}) to obtain a new
surface $\Sigma_1$ of genus 1 and an orientation-preserving
embedding $\psi_{1,1}\co S^1\times D^1\to \Sigma_1$. In fact we can
get a cobordism $(W_{0,1},\Sigma_0,\Sigma_1)$ with a
``gradient-like vector field'' $\xi$ for a Morse function
$f\co W_{0,1}\to [0,1]$. Here $f^{-1}(0) = \Sigma_0$, $f^{-1}(1) =
\Sigma_1$, and $f$ has a single critical point $p$ of index 1
with $f(p) = \frac{1}{2}$. We have that $\xi[f] >0$ away from $p$
and that in local coordinates near $p$, $f = \frac{1}{2} -{x_1}^2
+ {x_2}^2 + {x_3}^2$ and $\xi = -x_1\frac{\partial}{\partial
x_1}+x_2\frac{\partial}{\partial x_2}
+x_3\frac{\partial}{\partial x_3} $. The downward flow of $\xi$
from $p$ intersects $\Sigma_0$ in $\psi_{0,0}(S^0\times 0)$ and
the upward flow intersects $\Sigma_1$ in $\psi_{1,1}(S^1\times
0)$.

Choose an embedding $\psi_{0,1}\co S^0\times D^2\to \Sigma_1$ whose 
image 
is disjoint from $\psi_{1,1}(S^1\times D^1)$. Then we can repeat the 
process above to get another cobordism $(W_{1,2}, \Sigma_1, 
\Sigma_2)$ 
with Morse function $f\co  W_{1,2}\to [1,2]$ having a single critical 
point of index 1 at level $\frac{3}{2}$, and gradient-like vector 
field $\xi$ as before. 

Continuing in this way, we get a sequence of cobordisms 
$(W_{g,g+1},\Sigma_g,\Sigma_{g+1})$ between surfaces of genus 
difference 1, with Morse functions $f\co W_{g,g+1}\to [g,g+1]$ 
and gradient-like vector fields $\xi$. To each $\Sigma_g$, $g\geq 1$, is 
also associated a pair of embeddings $\psi_{i,g}\co S^i\times D^{2-i}\to 
\Sigma_g$, $i=0,1$. These embeddings have disjoint images, and are 
orientation-preserving with respect to the given, fixed orientations 
on the $\Sigma_g$. Note that the orientation on $\Sigma_g$ induced 
by $W_{g,g+1}$ is opposite to the given one, so the map 
$\psi_{0,g}\co S^0\times D^2\to -\Sigma_g=\partial_- W_{g,g+1}$ is
orientation-reversing.

Since the surfaces $\Sigma_g$ 
are all standard, we have a natural way to compose $W_{g-1,g}$ and 
$W_{g,g+1}$ to produce a cobordism $W_{g-1,g+1} = W_{g-1,g} + 
W_{g,g+1}$
 with a Morse function 
to $[g-1,g+1]$ having two index-1 critical points. Furthermore, by 
replacing $f$ by $-f$ we can obtain cobordisms $(W_{g+1,g}, 
\Sigma_{g+1}, \Sigma_g)$ with Morse functions having a single 
critical point of index 2, and these cobordisms may be naturally 
composed with each other or with the original index-1 cobordisms 
obtained before (after appropriately adjusting the values of the 
corresponding Morse functions), whenever such composition makes sense. 
We may think of $W_{g+1,g}$ as being simply $W_{g,g+1}$ with the 
opposite orientation.

In particular, we can fix integers $g,N\geq 0$ and proceed as 
follows. 
Beginning with $\Sigma_{g+N}$, compose the cobordisms 
$W_{g+N,g+N-1}, \ldots, W_{g+1, g}$ to form a ``standard'' 
compression body, and glue this with the 
composition $W_{g,g+1}+\cdots + W_{g+N-1,g+N}$ using the identity
map on $\Sigma_g$. The result is a cobordism $(W, \Sigma_{g+N},\Sigma_{g+N})$
and a Morse function $f\co  W\to \arr$ that we may rescale to have range $[-N,N]$,
having $N$ critical points each of index 1 and 2. By our
construction, the first half of this cobordism, $W_{g+N,g}$, is
identical with the second half, $W_{g,g+N}$: they differ only in
their choice of Morse function and associated gradient-like
vector field.

Now, by our construction the circles $\psi_{1,g+k}\co S^1\times 0\to 
f^{-1}(-k) = \Sigma_{g+k}\subset W$, $1\leq k\leq N$, all survive to
$\Sigma_{g+N}$ under downward flow of $\xi$. This is because the
images of $\psi_{1,q}$ and $\psi_{0,q}$ are disjoint for all $q$.
Thus on the ``lower'' copy of $\Sigma_{g+N}$ we have $N$ disjoint
primitive circles $c_1,\ldots, c_N$ that, under upward flow of
$\xi$, each converge to an index 2 critical point. Similarly,
(since $W_{g,g+N} = W_{g+N,g}$) the circles $\psi_{1,l}\co S^1\times
0 \to f^{-1}(k)=\Sigma_{g+k}\subset W$, $1\leq k\leq N$, survive to
$\Sigma_{g+N}$ under upward flow of $\xi$, and intersect the
``upper'' copy of $\Sigma_{g+N}$ in the circles $c_1,\ldots,c_N$.

Now suppose $h\co \partial_+W =\Sigma_{g+N}\to \Sigma_{g+N}=-\partial_-W$ 
is a diffeomorphism; 
then we can use $h$ to identify the boundaries $f^{-1}(-N)$, 
$f^{-1}(N)$ of $W$, and produce a manifold that we will denote by 
$M(g, N, h)$. Note that this manifold is entirely determined by the 
isotopy class of the map $h$, and that if $h$ preserves orientation 
then $M(g, N, h)$ is an orientable manifold having $b_1\geq 1$.

\begin{thm} Let $X$ be a closed oriented 3--manifold and $\phi\co X\to 
S^1$ a circle-valued Morse function with no critical points of index 
0 
or 3, and with $N$ critical points each of index 1 and 2. Assume 
that $[\phi]\in H^1(X;\zee)$ is of infinite order and indivisible. 
Arrange 
that $0<\arg\phi(p)<\pi$ for $p$ an index 1 critical point and 
$\pi<\arg\phi(q)<2\pi$ for $q$ an index 2 critical point, and let $S_g = 
\phi^{-1}(1)$, where $S_g$ has genus $g$. Then $X$ is
diffeomorphic to $M(g, N, h)$ for some
$h\co \Sigma_{g+N}\to\Sigma_{g+N}$ as above.
\label{stdthm}
\end{thm}

Note that $S_g$ has by construction the smallest genus among smooth
slices for $f$.

\proof By assumption $-1$ is a regular value of $\phi$, so 
$S_{g+N} = \phi^{-1}(-1)$ is a smooth orientable submanifold of $X$;
it is easy to see that $S_{g+N}$ is a closed surface of genus $g+N$.
Cut $X$ along $S_{g+N}$; then we obtain a cobordism $(W_\phi, S_-,
S_+)$ between two copies $S_\pm$ of $S_{g+N}$, and a Morse function $f\co 
W_\phi\to [-\pi,\pi]$ induced by $\arg\phi$. The critical points
of $f$ are exactly those of $\phi$ (with the same index), and by
our arrangement of critical points we have that $f(q)<0$ for any
index 2 critical point $q$ and $f(p)>0$ for any index 1 critical
point $p$. It is well-known that we can arrange for the critical
points of $f$ to have distinct values, and that in this case
$W_\phi$ is diffeomorphic to a composition of elementary
cobordisms, each containing a single critical point of $f$. For
convenience we rescale $f$ so that its image is the interval
$[-N,N]$ and the critical values of $f$ are the half-integers
between $-N$ and $N$. Orient each smooth level set $f^{-1}(x)$ by
declaring that a basis for the tangent space of $f^{-1}(x)$ is
positively oriented if a gradient-like vector field for $f$
followed by that basis is a positive basis for the tangent space
of $W_\phi$.

We will show that $W_\phi$ can be standardized by working ``from the 
middle out.'' Choose a gradient-like vector field $\xi_f$ for $f$, 
and 
consider $S_g = f^{-1}(0)$---the ``middle level'' of 
$W_\phi$, corresponding to $\phi^{-1}(1)$. There is exactly one 
critical point of $f$ in the region $f^{-1}([0,1])$, of index 1, and 
as above $\xi_f$ determines a ``characteristic embedding'' 
$\theta_{0,g}\co  S^0\times D^2\to S_g$. Choose a diffeomorphism 
$\Theta_0\co  S_g\to \Sigma_g$ such that $\Theta_0\circ\theta_{0,g} = 
\psi_{0,g}$; then it follows from \cite{milnor2}, Theorem 3.13, that 
$f^{-1}([0,1])$ is diffeomorphic to $W_{g,g+1}$ by some 
diffeomorphism $\Theta$
sending $\xi_f$ to $\xi$. (Recall that $\xi$ is the gradient-like
vector field fixed on $W_{g,g+1}$.)

Let $\Theta_1\co  S_{g+1}\to \Sigma_{g+1}$ be the restriction of  
$\Theta$
to $S_{g+1} = f^{-1}(1)$, and let $\mu_{0,g+1} = 
\Theta_1^{-1}\circ\psi_{0,g+1}\co S^0\times D^2\to S_{g+1}$.  Now 
$\xi_f$ induces an embedding $\theta_{0,g+1}\co S^0\times D^2\to 
S_{g+1}$, by considering downward flow from the critical point in 
$f^{-1}([1,2])$.  Since { any two orientation-preserving 
diffeomorphisms $D^2\to D^2$ are isotopic} and $S_{g+1}$ is connected, 
we have that $\mu_{0,g+1}$ and $\theta_{0,g+1}$ are isotopic.  It is 
now a simple matter to modify $\xi_f$ in the region 
$f^{-1}([1,1+\epsilon])$ using the isotopy, and arrange that 
$\theta_{0,g+1} = \mu_{0,g+1}$.  Equivalently, 
$\Theta\circ\theta_{0,g+1} = \psi_{0,g+1}$, so the theorem quoted 
above shows that $f^{-1}([1,2])$ is diffeomorphic to $W_{g+1, g+2}$.  
In fact, since the diffeomorphism sends $\xi_f$ to $\xi$, we get that 
$\Theta$ extends smoothly to a diffeomorphism $f^{-1}([0,2])\to 
W_{g,g+2}$.
 
Continuing in this way, we see that after successive modifications of 
$\xi_f$ in small neighborhoods of the levels $f^{-1}(k)$, $k = 
1,\ldots, N-1$, we obtain a diffeomorphism $\Theta\co  f^{-1}([0,N])\to 
W_{g,g+N}$ with $\Theta_*\xi_f = \xi$.
 
The procedure is entirely analogous when we turn to the ``lower 
half'' of $W_\phi$, but the picture is upside-down.  We have the 
diffeomorphism $\Theta_0\co  S_g\to \Sigma_g$, but before we can extend 
it to a diffeomorphism $\Theta\co  f^{-1}([-1,0])\to W_{g+1,g}$ we must 
again make sure the characteristic embeddings match.  That is, 
consider the map $\theta'_{0,g}\co  S^0\times D^2\to S_g$ induced by 
upward flow from the critical point, and compare it to 
$\Theta_0^{-1}\circ\psi_{0,g}$.  As before we can isotope $\xi_f$ in 
(an open subset whose closure is contained in) the region 
$f^{-1}([-\epsilon, 0])$ so that these embeddings agree, and we then 
get the desired extension of $\Theta$ to $f^{-1}([-1, N])$.  Then the 
procedure is just as before: alter $\xi_f$ at each step to make the 
characteristic embeddings agree, and extend $\Theta$ one critical 
point at a time.
 
Thus $\Theta\co  W_\phi \cong W = W_{g+N,g+N-1}+ \cdots+ 
W_{g+1,g}+W_{g,g+1}+\cdots + W_{g+N-1,g+N}$.  Since $W_\phi$ was 
obtained by cutting $X$, it comes with an identification $\iota\co  
S_+\to S_-$.  Hence $X\cong M(g, N, h)$ where $h = 
\Theta\circ\iota\circ\Theta^{-1}\co  \Sigma_{g+N}\to\Sigma_{g+N}$.
\endproof

\begin{rem} The identification $X\cong M(g, N, h)$ is not canonical, as it
depends on the initial choice of diffeomorphism $\phi^{-1}(1)\cong
\Sigma_g$, the final gradient-like vector field on
$W_\phi$ used to produce $\Theta$, as well as the function $\phi$.  As
with a Heegard decomposition, however, it is the existence of such a
structure that is important.
\end{rem}

\section{Preliminary calculations}

This section collects a few lemmata that we will use in the proof of
Theorem \ref{mainthm}.  Our main object here is to make the quantity
$[\zeta(F)\det(d)]_n$ a bit more explicit.

We work in the standardized setup of the previous section, identifying
$X$ with $M(g, N, h)$.  The motivation for doing so is mainly that our
invariants are purely algebraic---ie, homological---and the
standardized situation is very easy to deal with on this level.

Choose a metric $k$ on $X=M(g, N, h)$; then gradient flow with
respect to $k$ on $(W, \Sigma_{g+N},\Sigma_{g+N})$ determines
curves $\{c_i\}_{i=1}^N$ and $\{d_j\}_{j=1}^N$ on $\Sigma_{g+N}$,
namely $c_i$ is the intersection of the descending manifold of
the $i$th index-2 critical point with the lower copy of
$\Sigma_{g+N}$ and $d_j$ is the intersection of the ascending
manifold of the $j$th index-1 critical point with the upper copy
of $\Sigma_{g+N}$.

\begin{defn} The pair $(k,\phi)$ consisting of a metric $k$ on $X$ 
together with the Morse function $\phi\co X\to S^1$ is said to be 
{\em symmetric}
if the following conditions are satisfied.  Arrange the critical 
points of $\phi$ as in Theorem \ref{stdthm}, so that all critical 
points have distinct values. Write $W_\phi$ for the cobordism 
$X\setminus\phi^{-1}(-1)$, and $f\co W_\phi\to[-N,N]$ for the 
(rescaled) Morse function induced by $\phi$ as in the proof of Theorem 
\ref{stdthm}.  Write $I$ for the (orientation-reversing) involution 
obtained by swapping the factors in the expression $W_\phi \cong 
W_{g+N,g}\cup W_{g,g+N}$.  We require:
\begin{enumerate}
\item $I^*f = -f$.
\item For every $x\in W_{g+N,g}$ we have $(\grad f)_{I(x)} = -I_*(\grad 
f)_x$.
\end{enumerate}
\end{defn}

Symmetric pairs $(k,\phi)$ always exist: choose any metric on
$X$, and then in the construction used in the proof of Theorem 
\ref{stdthm}, take our gradient-like vector field $\xi_f$ to be a 
multiple of the gradient of $f$ with respect to that metric.  It is a 
straightforward exercise to see that the isotopies of $\xi_f$ needed in 
that proof may be obtained by modifications of the metric.

We use the term ``symmetric'' here because the gradient flows of the 
Morse function $f$ on the portions $W_{g+N,g}$ and $W_{g,g+N}$ are 
mirror images of each other. We will also say that the flow of $\grad 
f$ or of $\grad\phi$ is symmetric in this case.

Suppose $M(g,N,h)$ is endowed with a symmetric pair, and consider the 
calculation of $\zeta(F)\tau(M^*)$ in this case. Recall that $F$ is 
the return map of the flow of $\grad\phi$ from $\Sigma_g$ to itself 
(though $F$ is only partially defined due to the existence of 
critical points). Because of the symmetry of the flow, it is easy to see
that:
\begin{enumerate}
\item[(I)] The fixed points of iterates $F^k$ are in 1--1 correspondence 
with fixed points of iterates $h^k$ of the gluing map in the 
construction of $W$, and the Lefschetz signs of the fixed points 
agree. Indeed, if $h$ is sufficiently generic, we can 
assume that the set of fixed points of $h^k$ for $1\leq k\leq n$ (an 
arbitrary, but fixed, $n$) occur away from the $d_j$ (which 
agree with the $c_i$ under the identification 
$I$ by symmetry).  \item[(II)] The $(i,j)$th entry of the matrix of 
$d_M\co  M^1\to M^2$ in the Morse complex is given by the series
\begin{equation}
\sum_{k\geq 1} \langle {h^k}^* c_i,c_j\rangle t^{k-1},
\label{dformula}
\end{equation}
where $\langle \cdot,\cdot\rangle$ denotes the cup product pairing on 
$H^1(\Sigma_{g+N},\zee)$ and we have identified the curves $c_i$ with 
the Poincar\'e duals of the homology classes they represent.
\end{enumerate}

We should remark that a symmetric pair is not {\it a priori} 
suitable for calculating the invariant $\zeta(F)\tau(M^*)$ of 
Hutchings and Lee, since it is not generic. Indeed, for a
symmetric flow each index-2 critical point has a pair of upward
gradient flow lines into an index-1 critical point. However, this
is the only reason the flow is not generic: our plan now is to
perturb a symmetric metric to one which does not induce the
behavior of the flow just mentioned; then suitable genericity of
$h$ guarantees that the flow is Morse--Smale.

\begin{lemma} Assume that there are no ``short'' gradient flow lines 
between critical points, that is, every flow line between critical 
points intersects $\Sigma_g$ at least once.
Given a symmetric pair $(g_0,\phi)$ on $M(g, N, h)$ and suitable genericity 
hypotheses on $h$, there exists a $C^0$--small 
perturbation of $g_0$ to a metric $\tilde{g}$ such that for given $n\geq 0$ 
\begin{enumerate}
\item The gradient flow of $\phi$ with respect to $\tilde{g}$ is 
Morse--Smale; in particular the hypotheses of Theorem \ref{HLthm}
are satisfied. 
\item The quantity $[\zeta(F)\tau(M^*)]_m$, $m\leq
n$ does not change under this perturbation.
\end{enumerate}
\label{perturblemma}
\end{lemma}

We defer the proof of this result to Section \ref{lemmapfsec}.

\begin{rem} We can always arrange that there are no short gradient 
flow lines, at the expense of increasing $g = \mathrm{genus}(\Sigma_g)$. 
To see this, begin 
with $X$ and $\phi\co  X\to S^1$ as before, with $\Sigma_g =
\phi^{-1}(1)$ and the critical points arranged according to
index. Every gradient flow line then intersects $\Sigma_{g+N}$.
Now rearrange the critical points by an isotopy of $\phi$ that is
constant near $\Sigma_{g+N}$ so that the index-1 points occur in
the region $\phi^{-1}(\{e^{i\theta}|\pi<\theta<2\pi\})$ and the
index-2 points in the complementary region. This involves moving
all $2N$ of the critical points past $\Sigma_g$, and therefore
increasing the genus of the slice $\phi^{-1}(1)$ to $g+2N$; we
still have that every gradient flow line between critical points
intersects $\Sigma_{g+N}$.  Cutting $X$ along this new $\phi^{-1}(1)$ gives 
a cobordism $\tilde{W}$ between two copies of $\Sigma_{g+2N}$ and thus 
standardizes $X$ in the way we need while ensuring that there are no 
short flows.
\end{rem}

\begin{cor} The coefficients of the torsion $\tau(X,\phi)$ may be
calculated homologically, as the coefficients of the quantity
$\zeta(h)\tau(M^*_0)$ where $M^*_0$ is the Morse complex coming
from a symmetric flow.
\label{symcor}
\end{cor}

That is, we can use properties I and II of symmetric pairs to
calculate each coefficient of the right-hand side of (\ref{HLeqn}).

\begin{lemma}\label{detlemma}
If the flow of $\grad\phi$ is symmetric, the torsion 
$\tau(M^*)$ is represented by a polynomial whose $k$th coefficient is given by
\[
[\tau(M^*)]_k = \sum_{s_1+\cdots + s_N = k \atop \sigma\in{\mathfrak 
S}_N} (-1)^{\sgn(\sigma)}\langle {h^{s_1}}^*c_1,c_{\sigma(1)}\rangle\cdots 
\langle {h^{s_N}}^*c_N,c_{\sigma(N)}\rangle.
\]
\end{lemma}

\proof Since there are only two nonzero terms in the Morse 
complex, the torsion is represented by the determinant of the differential 
$d_M\co M^1\to M^2$.  Our task is to calculate a single coefficient of 
the determinant of this matrix of polynomials.  It will be convenient 
to multiply the matrix of $d_M$ by $t$; this multiplies $\det(d_M)$ by 
$t^N$, but $t^N\det(d_M)$ is still a representative for $\tau(M^*)$.  
Multiplying formula (\ref{dformula}) by $t$ shows
\begin{eqnarray*}
t^N\det(d_M) &=& \sum_{\sigma\in{\mathfrak S}_N} (-1)^{\sgn(\sigma)} 
\prod_i \left( \sum_k \langle {h^{k}}^* c_i,c_{\sigma(i)}\rangle 
t^k\right) \\
&=& \sum_{k} \sum_{\sigma\in{\mathfrak S}_N} \sum_{s_1+\cdots s_N = 
k}(-1)^{\sgn(\sigma)} \left( \prod_i \langle {h^{s_i}}^*c_i,
c_{\sigma(i)}\rangle\right) t^k
\end{eqnarray*}
and the result follows.\endproof

\section{Proof of Theorem \ref{mainthm}}

We are now in a position to explicitly calculate $\Tr \,\kappa_n$ using Theorem
\ref{wedgethm} and as a result prove Theorem \ref{mainthm}, assuming
throughout that $X$ is identified with $M(g, N, h)$ and the flow
of $\grad\phi$ is symmetric. Indeed, fix the nonnegative integer $n$ as
in Section \ref{tqftsec} and consider the cobordism $W_\phi$ as
above, identified with a composition of standard elementary
cobordisms. Using Theorem \ref{wedgethm} we see that the first
half of the cobordism, $W_{g+N,g}=f^{-1}([0,N])$, induces the map:
\begin{eqnarray*}
A_1\co  H^*(\Sigma_{g+N}^{(n+N)})&\to&
H^*(\Sigma_{g}^{(n)})\\
\alpha&\mapsto&\iota_{c_N^*}\cdots\iota_{c_1^*}\,\alpha
\end{eqnarray*}
The second half, $f^{-1}([N,2N])$, induces:
\begin{eqnarray*}
A_2\co  H^*(\Sigma_g^{(n)})&\to& H^*(\Sigma_{g+N}^{(n+N)})\\
\beta&\mapsto& c_1^*\wedge\cdots\wedge c_N^*\wedge\beta
\end{eqnarray*}
To obtain the map $\kappa_n$ we compose the above with the gluing
map $h^*$ acting on the symmetric power $\Sigma_{g+N}^{(n+N)}$. The
alternating trace $\Tr \,\kappa_n$ is then given by $\Tr(h^*\circ
A_2\circ A_1)$.

Following MacDonald \cite{MD}, we can take a monomial basis for
$H^*(\Sigma_g^{(n)})$.  Explicitly, if $\{x_i\}_{i = 1}^{2g}$ is a
symplectic basis for $H^1(\Sigma)$ having $x_i\cup x_{j+g} =
\delta_{ij}$ for $1\leq i,j\leq g$, and $x_i\cup x_j = 0$ for other
values of $i$ and $j$, $1\leq i<j\leq 2g$, and $y$ denotes the
generator of $H^2(\Sigma_g)$ coming from the orientation class, the
expression (\ref{symprodcohom}) shows that the set
\[
B_g^{(n)} = \{\alpha\} = \{x_Iy^q = x_{i_1}\wedge\cdots\wedge x_{i_k}\cdot
y^q| I = \{i_1<\ldots<i_k\}\subset\{1,\ldots,2g\}\},
\]
where $q=1,\ldots,n$ and $k=0,\ldots,n-q$, forms a basis for
$H^*(\Sigma_g^{(n)})$.  We take $H^*(\Sigma_{g+k}^{(n+k)})$ to have
similar bases $B_{g+k}^{(n+k)}$, using the images of the $x_i$ under the inclusion
$i\co H^1(\Sigma_{g+k-1})\to H^1(\Sigma_{g+k})$ constructed in section
\ref{tqftsec}, the (Poincar\'e duals of the) curves $c_1,\ldots,c_k$, 
and (the Poincar\'e duals of) some chosen dual curves $d_i$ to the $c_i$ as a 
basis for $H^1(\Sigma_{g+k})$.  Our convention is that $c_i\cup d_j = 
\delta_{ij}$, where we now identify $c_i$, $d_j$ with their Poincar\'e duals.

The dual basis for $B_{g+k}^{(n+k)}$ under the cup product pairing will be denoted
${B}_{n+k}^{\circ} = \{\alpha^{\circ}\}$.  Thus $\alpha^\circ\cup \beta =
\delta_{\alpha\beta}$ for basis elements $\alpha$ and
$\beta$. By abuse of notation, we will write $B_g^{(n)}\subset
B_{h}^{(m)}$ for $g\leq h$ and $n\leq m$; this makes use of the
inclusions on $H^1(\Sigma_g)$ induced by our standard cobordisms.

With these conventions, we can write:
\begin{eqnarray*}
\Tr\,\kappa_n &=& \sum_{\alpha\in B_{g+N}^{(n+N)}}(-1)^{\deg(\alpha)}
\alpha^\circ\cup h^*\circ A_2\circ A_1(\alpha)\\
&=& \sum_{\alpha\in B_{g+N}^{(n+N)}} (-1)^{\deg(\alpha)}\alpha^\circ\cup
h^*(c_1\wedge\cdots\wedge c_N \, \iota_{c_N}\cdots\iota_{c_1}
\alpha)
\end{eqnarray*}
For a term in this sum to be nonzero, $\alpha$ must be of a particular
form.  Namely, we must be able to write $\alpha = d_1\wedge\cdots
\wedge d_N\wedge \beta$ for some $\beta\in B_g^{(n)}$.  The sum then can
be written:
\begin{eqnarray}
&=& \sum_{\beta\in B_g^{(n)}}
(-1)^{\deg(\beta)+N}(d_1\wedge\cdots \wedge d_N\wedge \beta)^\circ \cup
h^*(c_1\wedge\cdots \wedge c_N\wedge \beta)
\label{traceform}
\end{eqnarray}
In words, this expression asks us to find the coefficient of
$d_1\wedge \cdots \wedge d_N\wedge \beta$ in the basis expression of
$h^*(c_1\wedge \cdots \wedge c_N \wedge\beta)$, and add up the results
with particular signs.  Our task is to express this coefficient in
terms of intersection data among the $c_i$ and the Lefschetz numbers
of $h$ acting on the various symmetric powers of $\Sigma_g$. 

Consider the term of (\ref{traceform}) corresponding to $\beta = 
x_Iy^q$ for $I=\{i_1,...,i_k\}\subset\{1,...,2g\}$ and $x_I = 
x_{i_1}\wedge\cdots\wedge x_{i_k}$. The coefficient of 
$d_1\wedge\cdots\wedge d_N\wedge x_Iy^q$ in the basis expression of 
$h^*(c_1\wedge\cdots\wedge c_N\wedge x_Iy^q)$ is computed by pairing 
each of $\{c_1,...,c_N,x_{i_1},...,x_{i_k}\}$ with each of 
$\{d_1,...,d_N,x_{i_1},...,x_{i_k}\}$ in every possible way, and 
summing the results with signs corresponding to the permutation 
involved. To make the notation a bit more compact, for given $I$ let $\bar{I} = 
\{1,...,N,i_1,...,i_k\}$ and write the elements of $\bar{I}$ as 
$\{\bar{\imath}_m\}_{m=1}^{N+k}$. Likewise, set $\bar{I}' = 
\{N+1,...,2N,i_1,...,i_k\} = 
\{\bar{\imath}'_1,...,\bar{\imath}'_{N+k}\}$.

Write $\{\xi_i\}_{i=1}^{2N+2g}$ for our basis of $H^1(\Sigma_{g+N})$: 
\begin{eqnarray*}
& \xi_1 = c_1, \quad\cdots,\quad \xi_N = c_N,\,\,\xi_{N+1} = 
d_1,\quad\cdots,\quad \xi_{2N} = d_N &\\
&\xi_{2N+1} =x_1,\quad\cdots,\quad \xi_{2N+2g} = x_{2g}&
\end{eqnarray*}
and let $\{\xi_i'\}$ be the dual basis: $\langle \xi_i,\xi'_j 
\rangle = \delta_{ij}$. Define $\zeta_i = h^*(\xi_i)$. 

Then since $\deg\beta = |I| = k$  modulo 2, the term of 
(\ref{traceform}) corresponding to $\beta = x_Iy^q$ is 
\begin{equation}
(-1)^{k+N} \sum_{\sigma\in\mathfrak{S}_{k+N}} (-1)^{\sgn(\sigma)} 
\langle \zeta_{\bar{\imath}_1}, 
\xi'_{\bar{\imath}_{\sigma(1)}'}\rangle \cdots \langle 
\zeta_{\bar{\imath}_{k+N}}, \xi'_{\bar{\imath}'_{\sigma(k+N)}}\rangle,
\label{coeffform}
\end{equation}
and (\ref{traceform}) becomes
\begin{equation}
\Tr\,\kappa_n =\sum_{k=0}^{\mathrm{min}(n,2g+2N)}\hspace*{-1em}
(2(n-k)+1)\hspace{-2em}\sum_{I\subset\{2N+1,\ldots,2N+2g\}\atop
|I| = k}\hspace*{-1.5em} [\mbox{formula (\ref{coeffform})}].
\label{traceform2}
\end{equation}

Here we are using the fact that for each $k =
0,\ldots,\mathrm{min}(n,2g+2N)$ the space
$\Lambda^kH^1(\Sigma_{g+N})$ appears in $H^*(\Sigma^{(n)})$
precisely $2(n-k)+1$ times, each in cohomology groups of all the
same parity.

Note that from (\ref{traceform}) we can see that the result is unchanged
if we allow not just sets $I\subset\{2N+1,\ldots,2N+2g\}$ in our
sum as above, but extend the sum to include sets $I =
\{i_1,\ldots,i_k\}$, where $i_1\leq\cdots \leq i_k$ and each
$i_j\in\{1,\ldots,2N+2g\}$. That is, we can allow $I$ to include
indices referring to the $c_i$ or $d_i$, and allow repeats: terms
corresponding to such $I$ contribute 0 to the sum. Likewise, we
may assume that the sum in (\ref{traceform2}) is over $k =
0,\ldots,n$ since values of $k$ larger than $2g+2N$ necessarily
involve repetitions in $\bar{I}$.

Consider the permutations $\sigma\in\mathfrak{S}_{k+N}$ used in the 
above. The fact that the first $N$ elements of $\bar{I}$ and 
$\bar{I}'$ are distinguished (corresponding to the $c_j$ and $d_j$, 
respectively) gives such permutations an additional structure. 
Indeed, writing $A =\{1,...,N\}\subset\{1,...,N+k\}$, let $\bar{A}$ 
denote the orbit of $A$ under powers of $\sigma$, and set $B = 
\{1,...,N+k\}\setminus \bar{A}$. Then $\sigma$ factors into a 
product $\sigma = \rho\cdot\tau$ where $\rho = \sigma|_{\bar{A}}$ 
and $\tau = \sigma|_B$. By construction, $\rho$ has the property that 
the orbit of $A$ under $\rho$ is all of $\bar{A}$. Given any 
integers $0\leq m\leq M$, we let $\mathfrak{S}_{M;m}$ denote the collection 
of permutations $\alpha$ of $\{1,...,M\}$ such that the orbit of 
$\{1,...,m\}$ under powers of $\alpha$ is all of $\{1,...,M\}$.  The 
discussion above can be summarized by saying that if $\bar{A} = 
\{a_1,...,a_N,a_{N+1},...,a_{N+r}\}$ (where $a_i = i$ for $i=1,...,N$) 
and $B = \{b_1,...,b_t\}$ then $\sigma$ preserves each of $\bar{A}$ 
and $B$, and $\sigma(\bar{A})= \{a_{\rho(1)},...,a_{\rho(N+r)}\}$, 
$\sigma(B) = \{b_{\tau(1)},...,b_{\tau(t)}\}$ for some 
$\rho\in\mathfrak{S}_{N+r;N}$, $\tau\in\mathfrak{S}_t$.  Furthermore, 
$\sgn(\sigma) = \sgn(\rho)+ \sgn(\tau)$ mod 2.

Finally, for $\rho\in \mathfrak{S}_{N+r;N}$ as above, we define
\[
s_i = \min\{m>0|\rho^m(i)\in\{1,...,N\}\}.
\]
The definition of $\mathfrak{S}_{N+r;N}$ implies that $\sum_{i=1}^N 
s_i = r+N$.

In (\ref{traceform2}) we are asked to sum over all sets $I$ with $|I|=k$ and all 
permutations $\sigma\in\mathfrak{S}_{N+k}$ of the subscripts of 
$\bar{I}$ and $\bar{I}'$. From the preceding remarks, this is
equivalent to taking a sum over all sets
$\bar{A}\supset\{1,...,N\}$ and $B$ with $|\bar{A}|+|B|=N+k$, and
all permutations $\rho$ and $\tau$,
$\rho\in\mathfrak{S}_{N+r;N}$, $\tau\in\mathfrak{S}_t$ (where
$|\bar{A}| = N+r$, $|B| = t$). Since we are to sum over all $I$
and $k$ and allow repetitions, we may replace $\bar{I}$ by
$\bar{A}\cup B$, meaning we take the sum over all $\bar{A}$ and
$B$ and all $\rho$ and $\tau$ as above, and eliminate reference
to $I$. Thus, we replace $\xi_{\bar{\imath}_{a_j}}$ by
$\xi_{a_j}$ and $\xi_{\bar{\imath}_{a_j}'}$ by $\xi_{a_j'}$ if we
define $\bar{A}' =
\{N+1,...,2N\}\cup(\bar{A}\setminus\{1,...,N\})$. (Put another
way, pairs $(\bar{I},\sigma)$ are in 1--1 correspondence with
4--tuples $(\bar{A}, B, \rho, \tau)$.) Then we can write
$\Tr\,\kappa_n$ as:
\begin{eqnarray*}
&&\hspace*{-2em}\sum_{k=0}^n
(2(n-k)+1)(-1)^{k+N}\hspace{-1em} \sum_{\bar{A},B\atop |\bar{A}|+|B| =
k+N}
\sum_{\rho\in\mathfrak{S}_{|A|;N}\atop\tau\in\mathfrak{S}_{|B|}}
(-1)^{\sgn(\rho)}\hspace{-1.1em}\prod_{i=1,\ldots,N\atop
m=0,\ldots,s_i-1}\langle
\zeta_{a_{\rho^m(i)}},\xi'_{a'_{\rho^{m+1}(i)}}\rangle\\
&&\hspace*{2in}\times
(-1)^{\sgn(\tau)}\prod_{r=1}^{|B|}\langle\zeta_{b_r},\xi'_{b_{\tau(r)}}\rangle
\end{eqnarray*}
Carrying out the sum over all $B$ of a given size $t$ and all
permutations $\tau$, this becomes:
\begin{eqnarray*}
&&\hspace*{-2em}\sum_{k=0}^n \sum_{\bar{A};
|\bar{A}| = k+N-t\atop t=0,\ldots,k}
\sum_{\rho\in\mathfrak{S}_{|A|;N}}
(-1)^{\sgn(\rho)+k+N}(2(n-k)+1)\hspace{-1em}\prod_{i=1,\ldots,N\atop
m=0,\ldots,s_i-1}\langle
\zeta_{a_{\rho^m(i)}},\xi'_{a'_{\rho^{m+1}(i)}}\rangle\\
&&\hspace*{2in}\times
\mathrm{tr}(h^*|_{\Lambda^tH^1(\Sigma_{g+N})})
\end{eqnarray*}
Reordering the summations so that the sum over $\bar{A}$ is on
the outside and the sum on $t$ is next, we find that $k =
|\bar{A}|-N+t$ and the expression becomes:
\begin{eqnarray*}
&&\hspace*{-2em} \sum_{\bar{A}\atop |\bar{A}|-N =
0,\ldots,n}\sum_{t=0}^{n-(|\bar{A}|-N)}
(-1)^{|\bar{A}|+\sgn(\rho)}\sum_{\rho\in\mathfrak{S}_{|A|;N}}
\prod_{i=1,\ldots,N\atop m=0,\ldots,s_i-1}\langle
\zeta_{a_{\rho^m(i)}},\xi'_{a'_{\rho^{m+1}(i)}}\rangle\\
&&\hspace*{1in}\times (-1)^t(2[n-(t-(|\bar{A}|-N))]+1)
\mathrm{tr}(h^*|_{\Lambda^tH^1(\Sigma_{g+N})})
\end{eqnarray*}
Again using the fact that $\Lambda^tH^1(\Sigma_{g+N})$ appears exactly
$2(|\bar{A}|-t)+1$ times in $H^*(\Sigma^{(|\bar{A}|-N)})$ and 
writing $|\bar{A}| = N+r$, we can
carry out the sum over $t$ to get that $\Tr\,\kappa_n$ is:
\[
\sum_{r=0}^n \left[\sum_{\bar{A}\atop|\bar{A}|-N =r}
\sum_{\rho\in\mathfrak{S}_{r+N;N}}
(-1)^{\sgn(\rho)+|\bar{A}|}\hspace{-1em}\prod_{i=1,\ldots,N \atop
m=0,\ldots,s_i-1}
\langle\zeta_{a_{\rho^m(i)}},\xi'_{a'_{\rho^{m+1}(i)}}\rangle \right]\cdot
L(h^{(n-r)})
\]
Here $L(h^{(n-r)})$ is the Lefschetz number of $h$ acting on
the $(n-r)$th symmetric power of $\Sigma_{g+N}$ which, as remarked 
in (\ref{zetasym}), is the $(n-r)$th coefficient of $\zeta(h)$. In 
view of Corollary \ref{symcor}, we will be done if we show
that the quantity in brackets is the $r$th coefficient of the 
representative $t^N\det(d_M)$ of $\tau(M^*)$.  Recalling the definition 
of $\bar{A}$, $\zeta_i$, and $\xi_i$, note that the terms that we are 
summing in the brackets above are products over all $i$ of formulae 
that look like
\begin{equation}
\langle c_i,\xi'_{a'_{\rho(i)}}\rangle\langle
h^*(\xi_{a_{\rho(i)}}),\xi'_{a'_{\rho^2(i)}}\rangle\cdots\langle
h^*(\xi_{\rho^{s_i-1}(i)}),c_{\tilde{\rho}(i)}\rangle
\label{example}
\end{equation}
where $\tilde{\rho}(i)\in\{1,\ldots,N\}$ is defined to be
$\rho^{s_i}(i)$. If we sum this quantity over all $\bar{A}$ and all
$\rho$ that induce the same permutation $\tilde{\rho}$ of
$\{1,\ldots,N\}$, we find that (\ref{example}) becomes simply $\langle
{h^*}^{s_i}(c_i),c_{\tilde{\rho}(i)}\rangle$. Therefore the
quantity in brackets is a sum of terms like
\[
(-1)^{\sgn(\rho)+r+N}\langle{h^*}^{s_1}c_1,c_{\tilde{\rho}(1)}\rangle\cdots\langle
{h^*}^{s_N}(c_N),c_{\tilde{\rho}(N)}\rangle,
\]
where we have fixed $s_1,\ldots,s_N$ and $\tilde{\rho}$ and carried out the
sum over all $\rho$ such that
\begin{enumerate}
\item $\mathrm{min}\{m>0|\rho^m(i)\in\{1,\ldots,N\}\} = s_i$, and
\item The permutation $i\mapsto \rho^{s_i}(i)$ of
$\{1,\ldots,N\}$ is $\tilde{\rho}$.
\end{enumerate}
(As we will see, $\sgn(\rho)$ depends only on
$\tilde{\rho}$ and $|\bar{A}|$.) It remains to sum over partitions
$s_1+\cdots +s_N$ of $s = |\bar{A}| = r+N$ and over permutations
$\tilde{\rho}$. But from Corollary \ref{symcor} and Lemma
\ref{detlemma}, the result of those two summations is precisely
$[\tau(M^*)]_r$, if we can see just that $\sgn(\tilde{\rho}) =
\sgn(\rho) + |\bar{A}|$ mod 2. That is the content of the next lemma.

\begin{lemma}
Let $A = \{1,\ldots,N\}$ and $\bar{A} = \{1,\ldots,s\}$ for some
$s\geq N$. Let $\rho\in\mathfrak{S}_{s;N}$ and define
\[
\tilde{\rho}(i)\in\mathfrak{S}_N, \,\,\tilde{\rho}(i) = \rho^{s_i}(i)
\]
where $s_i$ is defined as above.
Then $\sgn(\rho) = \sgn(\tilde{\rho}) + m$ modulo 2.
\end{lemma}

\proof Suppose $\rho = \rho_1\cdots\rho_p$ is an expression
of $\rho$ as a product of disjoint cycles; we may assume that the
initial elements $a_1,\ldots,a_p$ of $\rho_1,\ldots,\rho_p$ are
elements of $A$ since $\rho\in\mathfrak{S}_{m;N}$. For convenience we 
include any 1--cycles among the $\rho_i$, noting that the only elements 
of $\bar{A}$ that may be fixed under $\rho$ are in $A$. It is easy to
see that cycles in $\rho$ are in 1--1 correspondence with cycles
of $\tilde{\rho}$, so the expression of $\tilde{\rho}$ as a
product of disjoint cycles is $\tilde{\rho} =
\tilde{\rho}_1\cdots\tilde{\rho}_p$ where each $\tilde{\rho}_i$
has $a_i$ as its initial element. For $a\in A$, define
\begin{eqnarray*}
n(a) &=& \min\{m>0| \rho^m(a)\in A\} \\
\tilde{n}(a) &=& \min\{m>0 | \tilde{\rho}^m(a) = a\}.
\end{eqnarray*}
Note that $n(a_i) = s_i$ for $i=1,...,N$, $\sum s_i = s$, 
and $\tilde{n}(a_i)$ is the length of the cycle
$\tilde{\rho}_i$. The cycles $\rho_i$ are of the form
\[
\rho_i =
(a_i\cdots\tilde{\rho}(a_i)\cdots\tilde{\rho}^2(a_i)\cdots\cdots
\tilde{\rho}^{\tilde{n}(a_i) - 1}(a_i) \cdots)
\]
where ``$\cdots$'' stands for some number of elements of
$\bar{A}$. Hence the cycles $\rho_i$ have length
\[
l(\rho_i)\hspace{1ex} = \sum_{m=0}^{\tilde{n}(a_i) - 1}
(n(\tilde{\rho}^m(a_i))+1) \hspace{1ex}=\hspace{1ex}
\tilde{n}(a_i) +\hspace{-.5em}
\sum_{m=0}^{\tilde{n}(a_i)-1}\hspace{-.5em}n(\tilde{\rho}^m(a_i)).
\]
Modulo 2, then, we have
\begin{eqnarray*}
\sgn(\rho) &=& \sum_{i=1}^p (l(\rho_i) - 1) \\
&=& \sum_{i=1}^p\left[\left(\tilde{n}(a_i) +
\sum_{m=0}^{\tilde{n}(a_i)-1}n(\tilde{\rho}^m(a_i))\right) - 1\right]
\\
&=& \sum_{i=1}^p (\tilde{n}(a_i)-1) +
\sum_{i=1}^p\sum_{m=0}^{\tilde{n}(a_i)-1} n(\tilde{\rho}^m(a_i)) \\
&=& \sgn(\tilde{\rho}) + s,
\end{eqnarray*}
since because $\rho\in\mathfrak{S}_{s;N}$ we have 
$\sum_{i=1}^p\sum_{m=0}^{\tilde{n}(a_i)-1} n(\tilde{\rho}^m(a_i)) = 
\sum_{i=1}^N s_i = s$.\endproof

\section{Proof of Theorem \ref{intthm}}

The theorem of Hutchings and Lee quoted at the beginning of this
work can be seen as (or more precisely, the logarithmic derivative
of formula (\ref{HLeqn}) can be seen as) a kind of Lefschetz
fixed-point theorem for partially-defined maps, specifically the
return map $F$, in which the torsion $\tau(M^*)$ appears as a
correction term (see \cite{HL1}). Now, the Lefschetz number of a
homeomorphism $h$ of a closed compact manifold $M$ is just the
intersection number of the graph of $h$ with the diagonal in
$M\times M$; such consideration motivates the proof of Theorem
\ref{HLthm} in \cite{HL1}. With the results of Section
\ref{stdsec}, we can give another construction.

Given $\phi\co X = M(g,N,h)\to S^1$ our circle-valued Morse function, cut along
$\phi^{-1}(-1)$ to obtain a cobordism $W_\phi$ between two copies
of $\Sigma_{g+N}$. Write $\gamma_i$, $i=1,\ldots,N$ for the
intersection of the ascending manifolds of the index-1 critical
points with $\partial_+W$ and $\delta_i$ for the intersection of
the descending manifolds of the index-2 critical points with
$\partial_-W$. Since the homology classes $[\gamma_i]$ and
$[\delta_i]$ are the same (identifying
$\partial_+W=\partial_-W=\Sigma_{g+N}$), we may perturb the
curves $\gamma_i$ and $\delta_i$ to be parallel, ie, so that they do not
intersect one another (or any other $\gamma_j$, $\delta_j$ for
$j\neq i$ either). Choose a complex structure on $\Sigma_{g+N}$
and use it to get a complex structure on the symmetric powers
$\Sigma_{g+N}^{(k)}$ for each $k$. Write $T_\gamma$ for the
$N$--torus $\gamma_1\times\cdots\times \gamma_N$ and let
$T_\delta = \delta_N\times\cdots\times\delta_1$. Define a function
\[
\psi\co  T_\gamma\times \Sigma_{g+N}^{(n)}\times T_\delta\to
\Sigma_{g+N}^{(n+N)}\times\Sigma_{g+N}^{(n+N)}
\]
by mapping the point $(q_1,\ldots,q_N,\sum p_i,q'_N,\ldots,q'_1)$ 
to $(\sum p_i + \sum q_j, \sum p_i + \sum q_j')$. 

The perhaps unusual-seeming orders on the $\delta_i$ and in the
domain of $\psi$ are chosen to obtain the correct sign in the sequel.

\begin{prop} $\psi$ is a smooth embedding, and $D =
\mathrm{Im}\psi$ is a totally real submanifold of
$\Sigma_{g+N}^{(n+N)}\times\Sigma_{g+N}^{(n+N)}$.
\end{prop}

The submanifold $D$ plays the role of the diagonal in the
Lefshetz theorem. 

\proof That $\psi$ is one-to-one is clear since the
$\gamma_i$ and $\delta_j$ are all disjoint. For smoothness, we
work locally. Recall that the symmetric power $\Sigma_g^{(k)}$ is
locally isomorphic to $\cee^{(k)}$, and a global chart on the latter is
obtained by mapping a point $\sum w_i$ to the coefficients of the
monic polynomial of degree $k$ having zeros at each $w_i$. Given a
point $(\sum p_i +\sum q_j, \sum p_i+\sum q_j')$ 
of $\mathrm{Im}(\psi)$ we can choose a coordinate chart on 
$\Sigma_{g+N}$ containing all the points $p_i,q_j,q_j'$ so that the $\gamma_i$ and 
$\delta_j$ are described by disjoint curves in $\cee$.  Thinking of 
$q_j\in\gamma_j\subset\cee\cong\cee^{(1)}$ and simlarly for $q_j'$, we 
have that locally $\psi$ is just the multiplication map:
\begin{eqnarray*}
&&\hspace*{-.5in}\left(
(z-q_1),\ldots,(z-q_N),\prod_{i=1}^n(z-p_i),(z-q_1'),\ldots,(z-q_N')\right)\\
&&\hspace*{.5in}\mapsto \left(\prod_{i=1}^n
(z-p_i)\prod_{j=1}^N(z-q_j),\hspace{1ex}\prod_{i=1}^n(z-p_i)\prod_{j=1}^N(z-q_j')\right)
\end{eqnarray*}
It is clear that the coefficients of the polynomials on the right hand side
depend smooth\-ly on the coefficients of the one on the left
and on the $q_j$, $q_j'$.

On the other hand, if $(f(z),g(z))$ are the polynomials whose
coefficients give the local coordinates for a point in
$\mathrm{Im}(\psi)$, we know that $f(z)$ and $g(z)$ share exactly
$n$ roots since the $\gamma_i$ and $\delta_j$ are disjoint. If
$p_1$ is one such shared root then we can write $f(z) =
(z-p_1)\tilde{f}(z)$ and similarly for $g(z)$, where
$\tilde{f}(z)$ is a monic polynomial of degree $n+N-1$ whose
coefficients depend smoothly (by polynomial long division!) on
$p_1$ and the coefficients of $f$. Continue factoring in this way
until $f(z) = f_0(z)\prod_{i=1}^n(z-p_i)$, using the fact that
$f$ and $g$ share $n$ roots to find the $p_i$. Then $f_0$ is a
degree $N$ polynomial having one root on each $\gamma_i$, hence
having all distinct roots.  Those roots (the $q_j$) therefore depend smoothly on 
the coefficients of $f_0$, which in turn depend smoothly on the 
coefficients of $f$.  Hence $D$ is smoothly embedded.

That $D$ is totally real is also a local calculation, and is a
fairly straightforward exercise from the definition.\endproof

We are now ready to prove the ``algebraic'' portion of Theorem
\ref{intthm}.

\begin{thm}\label{algintthm}
Let $\Gamma$ denote the graph of the map $h^{(n+N)}$
induced by the gluing map $h$ on the symmetric product
$\Sigma_{g+N}^{(n+N)}$. Then
\[
D.\Gamma = \Tr\,\kappa_n.
\]
\end{thm}

\proof Using the notation from the previous section, we
have that in cohomology the duals of $D$ and $\Gamma$ are
\begin{eqnarray*}
D^* &=& \sum_{\beta\in B_{g+N}^{(n)}} (-1)^{\epsilon_1(\beta)}
(c_1\wedge\cdots \wedge c_N\wedge \beta^{\circ})\times
(c_1\wedge\cdots \wedge c_N\wedge\beta) \\
\Gamma^*&=& \sum_{\alpha\in B_{g+N}^{(n+N)}}(-1)^{\deg(\alpha)}
\alpha^\circ\times {h^*}^{-1}(\alpha).
\end{eqnarray*}
 Here $\epsilon_1(\beta) =
\deg(\beta)(N+1)+\frac{1}{2}N(N-1)$. Indeed, since
the diagonal is the pushforward of the graph by $1\times
h^{-1}$, we get that the dual of the graph is the pullback of the
diagonal by $1\times h^{-1}$. We will find it convenient to write
\[
D^* = \sum_\beta (-1)^{\epsilon_1(\beta) + \epsilon_2(\beta)}
(c_1\wedge\cdots \wedge c_N\wedge \beta)\times (c_1\wedge\cdots
\wedge c_N\wedge\beta^{\circ}),
\]
by making the substitution $\beta\mapsto\beta^{\circ}$ in the
previous expression. Since $\beta^{\circ\circ} =
\pm\beta$, the result is still a sum over the monomial basis with
an additional sign denoted by $\epsilon_2$ in the above but which
we will not specify.

Therefore the
intersection number is
\begin{equation}\begin{split}
D^*\cup \Gamma^* =& \sum_{\alpha,\beta}
(-1)^{\epsilon_1 + \epsilon_2+ \epsilon_3(\alpha,\beta)}\\&
(\alpha^\circ\cup(c_1\wedge\cdots\wedge c_N\wedge\beta))\times
({h^*}^{-1}\alpha\cup(c_1\wedge\cdots\wedge
c_N\wedge\beta^\circ))\end{split}
\label{intformula1}
\end{equation}
where $\epsilon_3(\alpha,\beta) = \deg(\alpha)(1+\deg(\beta) +
N)$. Since this is a sum over a monomial basis $\alpha$, the
first factor in the cross product above vanishes unless $\alpha = 
c_1\wedge\cdots\wedge c_N\wedge\beta$, and in that case is 1.  
Therefore $\deg(\alpha) = \deg(\beta)+N$, which gives 
$\epsilon_3(\alpha,\beta) \equiv 0$ mod 2, and (\ref{intformula1}) becomes
\begin{eqnarray}
D^*\cup\Gamma^* &=& \sum_\beta (-1)^{\epsilon_1 + \epsilon_2}
{h^*}^{-1}(c_1\wedge\cdots\wedge c_N\wedge\beta)\cup
(c_1\wedge\cdots\wedge c_N\wedge\beta^\circ)\nonumber\\
&=& \sum_\beta (-1)^{\epsilon_1+\epsilon_2}
(c_1\wedge\cdots\wedge c_N\wedge\beta) \cup
h^*(c_1\wedge\cdots\wedge c_N\wedge\beta^\circ)\nonumber\\
&=& \sum_\beta (-1)^{\epsilon_1} (c_1\wedge\cdots\wedge
c_N\wedge\beta^\circ) \cup h^*(c_1\wedge\cdots\wedge c_N\wedge\beta)
\label{intformula2}
\end{eqnarray}
where we have again used the substitution
$\beta\mapsto\beta^\circ$ and therefore cancelled the sign $\epsilon_2$.
Now, some calculation using the cup product structure of
$H^*(\Sigma_{g+N}^{(n+N)})$ derived in \cite{MD} shows that
\[
c_1\wedge\cdots\wedge c_N\wedge\beta^\circ =
(-1)^{\epsilon_4(\beta)}(d_1\wedge\cdots\wedge d_N\wedge\beta)^\circ.
\]
where $\epsilon_4(\beta) = N\deg(\beta) +
\frac{1}{2}N(N+1) \equiv \epsilon_1(\beta) + \deg(\beta) + N$ mod 2.
Note that ${(\cdot)}^\circ$ refers to duality in
$H^*(\Sigma_{g+N}^{(n)})$ on the left hand side and in
$H^*(\Sigma_{g+N}^{(n+N)})$ on the right. Returning with this to
(\ref{intformula2}) gives
\[
D^*\cup\Gamma^* = \sum_\beta (-1)^{\deg(\beta)+N}
(d_1\wedge\cdots\wedge d_N\wedge\beta)^\circ\cup
h^*(c_1\wedge\cdots\wedge c_N\wedge\beta),
\]
which is $\Tr\,\kappa_n$ by (\ref{traceform}). Theorem
\ref{algintthm} follows.\endproof

To complete the proof of Theorem \ref{intthm}, we recall that we
have already shown that $D$ is a totally real submanifold of
$\Sigma_{g+N}^{(n+N)}\times\Sigma_{g+N}^{(n+N)}$. The graph of
$h^{(n+N)}$, however, is not even smooth unless $h$ is an
automorphism of the chosen complex structure of $\Sigma_{g+N}$:
in general the set-theoretic map induced on a symmetric power by
a diffeomorphism of a surface is only Lipschitz continuous.
Salamon \cite{S} has shown that if we choose a path of
complex structures on $\Sigma$ between the given one $J$ and
$h^*(J)$, we can construct a symplectomorphism of the moduli
space ${\cal M}(\Sigma,J)\cong\Sigma_{g+N}^{(n+N)}$ that is
homotopic to the induced map $h^{(n+N)}$. Hence $\Gamma$ is
homotopic to a Lagrangian submanifold of
$\Sigma_{g+N}^{(n+N)}\times-\Sigma_{g+N}^{(n+N)}$. Since
Lagrangians are in particular totally real, and since
intersection numbers do not change under homotopy, Theorem
\ref{intthm} is proved.

\section{Proof of Lemma \ref{perturblemma}}
\label{lemmapfsec}

We restate the lemma:

{\sl Assume that there are no ``short'' gradient flow lines 
between critical points, that is, every flow line between critical 
points intersects $\Sigma_g$ at least once.
Given a symmetric pair $(g_0,\phi)$ on $M(g, N, h)$ and suitable genericity 
hypotheses on $h$, there exists a $C^0$--small 
perturbation of $g_0$ to a metric $\tilde{g}$ such that for given $n\geq 0$ 
\begin{enumerate}
\item The gradient flow of $\phi$ with respect to $\tilde{g}$ is 
Morse--Smale; in particular the hypotheses of Theorem \ref{HLthm}
are satisfied. 
\item The quantity $[\zeta(F)\tau(M^*)]_m$, $m\leq
n$ does not change under this perturbation.
\end{enumerate}}

\proof
Alter $g_0$ in a small neighborhood of $\Sigma_g\subset M(g, N, h)$ as 
follows, working in a half-collar neighborhood of $\Sigma_g$ diffeomorphic to 
$\Sigma_g\times (-\epsilon, 0]$ using the flow of $\grad_{g_0}\phi$ to 
obtain the product structure on this neighborhood.

Let $p_1,\ldots,p_{2N}$ denote the 
points in which the ascending manifolds (under gradient flow of $f$ 
with respect to the symmetric metric $g_0$) of the index-2 critical points 
intersect $\Sigma_g$ in $W_\phi$. Since $g_0$ is symmetric, these 
points are the same as the points $q_1,\ldots,q_{2N}$ in which the 
descending manifolds of the index-1 critical points intersect 
$\Sigma_g$. Let ${\cal O}$ denote the union of all closed orbits of 
$\grad\phi$ (with respect to $g_0$) of degree no more than $n$, and 
all gradient flow lines connecting index-1 to index-2 critical points. We may 
assume that this is a finite set. Choose small disjoint coordinate disks 
$U_i$ around each $p_i$ such that $U_i\cap ({\cal O}\cap \Sigma_g) = 
\emptyset$. 

In $U_i\times (-\epsilon,0]$, we may suppose the Morse function $f$ is given by 
projection onto the second factor, $(u,t)\mapsto t$, and the metric is 
a product $g_0 = g_{\Sigma_g}\oplus (1)$. Let ${X}_i$ be a 
nonzero constant 
vector field in the coordinate patch $U_i$ and $\mu$ a cutoff function 
that is equal to $1$ near $p_i$ and zero off a small neighborhood of 
$p_i$ whose closure is in $U_i$.  Let $\nu(t)$ be a bump function that 
equals 1 near $t = \epsilon/2$ and vanishes near the ends of the 
interval $(-\epsilon,0]$.  Define the vector field $v$ in the set 
$U_i\times (-\epsilon, 0]$ by $v(u,t) = \grad_{g_0}\phi + \nu(t)\mu(u) 
X(u)$.  Now define the metric $g_{X_i}$ in $U_i\times (-\epsilon,0]$
 by declaring that $g_{X_i}$ agrees with $g_0$ 
on tangents to slices $U_i\times\{t\}$, but that $v$ is orthogonal to the 
slices.  Thus, with respect to $g_{X_i}$, the gradient $\grad\phi$ is 
given by a multiple of $v(u,t)$ rather than $\partial/\partial t$.

It is easy to see that repelacing $g_0$ by $g_{X_i}$ in $U_i\times 
(-\epsilon,0]$ for each $i = 1,\ldots,2N$ produces a metric $g_X$ for 
which upward gradient flow of $\phi$ on $W_\phi$ does not connect index-2 
critical points to index-1 critical points with ``short'' gradient 
flow lines.  Elimination of gradient flows of $\phi$ from index-2 to index-1 
points that intersect $\Sigma_{g+N}$ is easily arranged by small 
perturbation of $h$, as are transverse intersection of ascending and 
descending manifolds and nondegeneracy of fixed points of $h$ and its 
iterates.  Hence the new metric $g_X$ satisfies condition (1) of the 
Lemma.

For condition (2), we must verify that we have neither created nor 
destroyed either closed orbits of $\grad\phi$ or flows from index-1 
critical points to index-2 critical points. The fact that no such 
flow lines have been destroyed is assured by our choice of neighborhoods 
$U_i$. We now show that we can choose the  vector fields $X_i$ such 
that no fixed points of $F^k$ are created, for $1\leq k\leq n$.

Let $F_1\co  \Sigma_g\to \Sigma_{g+N} = \partial_+W_\phi$ be the map induced by gradient 
flow with respect to $g_0$, defined away from the $q_j$, and let 
$F_2\co  \Sigma_{g+N} = \partial_-W_\phi\to\Sigma_g$ be the similar map 
from the bottom of the cobordism, defined away from the $c_j$. Then 
the flow map $F$, with respect to $g_0$, is given by the composition 
$F = F_2\circ h\circ F_1$ where this is defined.  The return map with 
respect to the $g_X$--gradient, which we will write $\tilde{F}$, is 
given by $F$ away from the $U_i$ and by $F + cX$ in the coordinates on 
$U_i$ where $c$ is a nonnegative function on $U_i$ depending on $\mu$ 
and $\nu$, vanishing near $\partial U_i$.

Consider the graph $\Gamma_{F^k}\subset \Sigma_g\times\Sigma_g$. Since 
$F^k$ is 
not defined on all of $\Sigma_g$ the graph is not closed, nor is its 
closure a cycle since $F^k$ in general has no continuous extension to 
all of $\Sigma$. Indeed, the boundary of $\Gamma_{F^k}$ is given by a 
union of products of ``descending slices'' (ie, the intersection of 
a descending manifold of a critical point with $\Sigma_g$) with 
ascending slices. Restrict attention to the neighborhood 
$U$ of $p$, where for convenience $p$ denotes any of the 
$p_1,\ldots,p_{2N}$ above. We have chosen $U$ so that there are no fixed 
points of $F^k$ in this neighborhood, ie, the graph and the diagonal 
are disjoint over $U$. If there is an open set around $\Gamma_{F^k}\cap 
(U\times U)$ that misses the diagonal $\Delta\subset U\times 
U$, then any sufficiently small choice of $X$ will keep $\Gamma_{F^k}$ 
away from $\Delta$ and therefore produce no new closed orbits of the 
gradient flow.  However, it may be that $\partial \Gamma_{F^k}$ has points 
on $\Delta$.  Indeed, if $c\subset\partial_+W_\phi = \Sigma_{g+N}$ is 
the ascending slice of the critical point corresponding to $p=q$, 
suppose $h^k(c)\cap c\neq \emptyset$.  Then it is not hard to see that 
$(p,p)\in\partial\Gamma_{F^k}$, and this situation cannot 
be eliminated by genericity assumptions on $h$.  Essentially, $p$ is 
both an ascending slice and a descending slice, so 
$\partial\Gamma_{F^k}$ can contain both $\{p\}\times(\mathrm{asc.  
slice})$ and $(\mathrm{desc.  slice})\times\{p\}$, and ascending and 
descending slices can have $p$ as a boundary point.

Our perturbation of $F$ using $X$ amounts, over $U$, to a ``vertical'' 
isotopy of $\Gamma_{F^k}\subset U\times U$. The question of whether there is an $X$ 
that produces no new fixed points is that of whether there is a 
vertical direction to move $\Gamma_{F^k}$ that results in the 
``boundary-fixed'' points like $(p,p)$ described above remaining 
outside of $int(\Gamma_{F^k})$. The existence of such a direction is 
equivalent to the jump-discontinuity of $F^k$ at $p$. This argument is 
easy to make formal in the case $k=1$, and for $k>1$ the ideas are the 
same, with some additional bookkeeping. We leave the general argument 
to the reader.

Turn now to the question of whether any new flow lines between 
critical points are created. Let $D = (h\circ F_1)^{-1}(\bigcup c_i)$ denote 
the first time that the descending manifolds of the critical points 
intersect $\Sigma_g$, and let $A = F_2\circ h (\bigcup c_i)$ be the 
similar ascending slices. Then except for short flows, the flow lines 
between critical points are in 1--1 correspondence with intersections 
of $D$ and $F^k(A)$, for various $k\geq 0$. We must show that our 
perturbations do not introduce new intersections between these sets. 
It is obvious from our constructions that only $F^k(A)$ is affected by 
the perturbation, since only $F_2$ is modified. 

Since there are no short flows by assumption, there are no 
intersections of $h^{-1}(c_j)$ with $c_i$ for any $i$ and $j$. This means 
that $D$ consists of a collection of embedded circles in $\Sigma_g$, 
where in general it may have included arcs connecting various $q_i$. 
Hence, we can choose our neighborhoods $U_i$ small enough that 
$U_i\cap D = \emptyset $ for all $i$, and therefore the perturbed 
ascending slices $\tilde{F}^k(A)$ stay away from $D$. Hence no new 
flows between critical points are created.

This concludes the proof of Lemma \ref{perturblemma}.\endproof

\end{document}

%% file: gtoutput.tex
\def\ifplaintex{\expandafter\ifx\csname documentclass\endcsname\relax}

\ifplaintex 
\else
\fi

\expandafter\ifx\csname beginpicture\endcsname\relax
\expandafter\ifx\csname documentclass\endcsname\relax
\input pictex \else
\input prepictex \input pictex \input postpictex \fi\fi

\def\gt{{\mathsurround=0pt\it $\cal G\mskip-2mu$eometry \&\ 
$\cal T\!\!$opology}}        %  journal title in recommended style

\def\gtp{{\mathsurround=0pt\it $\cal G\mskip-2mu$eometry \&\ 
$\cal T\!\!$opology $\cal P\!$ublications}}  % GT publications

\def\lognumber#1{\def\thelognumber{#1}}
\def\volumenumber#1{\def\thevolumenumber{#1}}
\def\papernumber#1{\def\thepapernumber{#1}}
\def\volumeyear#1{\def\thevolumeyear{#1}}

\def\pagenumbers#1#2{\def\startpage{#1}\def\finishpage{#2}}
\def\published#1{\def\publishdate{#1}}
\def\proposed#1{\def\theproposer{#1}}
\def\seconded#1{\def\theseconders{#1}}
\def\received#1{\def\receiveddate{#1}}

\def\accepted#1{\def\accepteddate{#1}}
\def\asciititle#1{\def\theasciititle{#1}}

\long\def\asciiabstract#1{\long\def\theasciiabstract{#1}}
\def\asciikeywords#1{\def\theasciikeywords{#1}}

\let\\\par\let\thelognumber\relax
\let\thevolumenumber\relax\let\thepapernumber\relax
\let\thevolumeyear\relax\let\thesamplenumber\relax\let\startpage\relax
\let\finishpage\relax\let\publishdate\relax\let\receiveddate\relax
\let\reviseddate\relax\let\accepteddate\relax\let\theasciititle\relax
\let\theasciiauthors\relax
\let\theasciiabstract\relax\let\theasciikeywords\relax
\let\theasciiemail\relax\let\theshortauthors\relax\let\theshorttitle\relax

\long\def\maketitlep{   % start of definition of \maketitlep

\count0=\startpage

\gt\hfill      %   Journal title (top left) 
\beginpicture
\setcoordinatesystem units <0.33truein, 0.33truein> point at 2.2 0.9
\setplotsymbol ({$\cal G$})
\plotsymbolspacing=9truept
\circulararc 315 degrees from 0 1 center at 0 0
\setplotsymbol ({$\cal T$})
\circulararc 315 degrees from 1 -1 center at 1 0
\endpicture
\break
{\small\ifx\thesamplenumber\relax % sample?  
Volume \else Sample
\fi\thevolumenumber\ (\thevolumeyear)
\startpage--\finishpage\nl
Published: \publishdate}
\vglue 0.5truein plus 0.4fil minus 0.1truein

{\parskip=0pt\leftskip 0pt plus 1fil\def\\{\par\smallskip}{\ifplaintex\large
\else\Large\fi\bf\thetitle}\par\medskip}   

\vglue 0pt plus 0.1fil 

{\parskip=0pt\leftskip 0pt plus 1fil\def\\{\par}{\sc\theauthors}
\par\medskip}

\vglue 0pt plus 0.1fil 

{\small\parskip=0pt\let\newline\\
{\leftskip 0pt plus 1fil\def\\{\par}{\sl\theaddress}\par}
\expandafter\ifx\theemail\relax    % email address?
\relax\else\vglue 5pt plus 0.02fil minus 2pt\def\\{\stdspace{\rm 
and}\stdspace} 
\cl{Email:\stdspace\tt\theemail}\fi
\ifx\theurl\relax                  % URL given?
\relax\else\vglue 5pt plus 0.02fil minus 2pt\def\\{\stdspace{\rm 
and}\stdspace}
\cl{URL:\stdspace\tt\theurl}\fi\par}

\vglue 7pt plus 0.3fil minus 3pt

{\bf Abstract}
\vglue 5pt plus 0.1fil minus 2pt

\theabstract

\vglue 7pt plus 0.3fil minus 3pt

{\bf AMS Classification numbers}\quad Primary:\quad \theprimaryclass

Secondary:\quad \thesecondaryclass

\vglue 5pt plus 0.3fil minus 2pt

{\bf Keywords}\quad \thekeywords

\vglue 10pt plus 0.5fil minus 5pt

{\small  Proposed: \theproposer\hfill Received: \receiveddate\nl
Seconded: \theseconders\hfill 
\ifx\reviseddate\relax                         % paper revised?
Accepted: \accepteddate                        % no
\else
Revised: \reviseddate                          % yes
\fi}
\eject
}       %  end of definition of \maketitlep

\let\maketitlepage\maketitlep
\let\maketitle\maketitlepage

\font\phead=cmsl9 scaled 950
\font=cmsl9 scaled 1050
\font\pnum=cmbx10 scaled 913
\font\lnum=cmbx10 
\font\pfoot=cmsl9 scaled 950
\font=cmsl9 scaled 1050
\ifplaintex
\headline{\vbox to 0pt{\vskip -4.5mm\line{\small\phead\ifnum
\count0=\startpage ISSN 1364-0380 (on line)
1465-3060 (printed) \hfill {\pnum\folio}\else\ifodd\count0\def\\{ }% 
\ifx\theshorttitle\relax\thetitle\else\theshorttitle\fi\hfill{\pnum\folio}
\else\def\\{ and }{\pnum\folio}\hfill\ifx\theshortauthors\relax\theauthors
\else\theshortauthors\fi\fi\fi}\vss}}
\footline{\vbox to 0pt{\vglue 0mm\line{\small\pfoot\ifnum\count0=\startpage
\copyright\ \gtp\hfill\else
\gt, Volume \thevolumenumber\ (\thevolumeyear)\hfill\fi}\vss
}}
\else
\makeatletter
\def\@oddhead{{\small\lhead\ifnum\count0=\startpage ISSN 1364-0380 (on line)
1465-3060 (printed) \hfill {\lnum\number\count0}\else\ifodd\count0
\def\\{ }\ifx\theshorttitle\relax \thetitle \else\theshorttitle\fi\hfill
{\lnum\number\count0}\else\def\\{ and }{\lnum\number\count0}
\hfill\ifx\theshortauthors\relax 
\theauthors\else\theshortauthors\fi\fi\fi}}\def\@evenhead{\@oddhead}
\def\@oddfoot{\small\lfoot\ifnum\count0=\startpage\copyright\ \gtp\hfill\else
\gt, Volume \thevolumenumber\ (\thevolumeyear)\hfill\fi}
\def\@evenfoot{\@oddfoot}
\makeatother
\fi

\newwrite\gtoutfile
\long\gdef\makeheadfile{  %%% start of definition of \makeheadfile
{\def\\{, }\def\s{ }
\immediate\openout\gtoutfile head.xxx
\immediate\write\gtoutfile{To: math@arxiv.org}
\immediate\write\gtoutfile{Subject: put or rep NNNNN:pppp}
\immediate\write\gtoutfile{--text follows this line--}
\immediate\write\gtoutfile{Proxy-for: \ifx\theasciiauthors\relax
\theauthors\else\theasciiauthors\fi\s<\ifx\theasciiemail\relax\theemail\else\theasciiemail\fi>}
\immediate\write\gtoutfile{\noexpand\\}
\immediate\write\gtoutfile{Authors: \ifx\theasciiauthors\relax
\theauthors\else\theasciiauthors\fi}
{\def\\{ }\immediate\write\gtoutfile{Title: \ifx\theasciititle\relax
\thetitle\else\theasciititle\fi}}
\immediate\write\gtoutfile{Subj-class: GT or SG or MG etc}
\immediate\write\gtoutfile{MSC-class: \theprimaryclass\ifx\thesecondaryclass\relax\else, \thesecondaryclass\fi}
\immediate\write\gtoutfile{Journal-ref: Geom. Topol. \thevolumenumber
(\thevolumeyear) \startpage-\finishpage}
\immediate\write\gtoutfile{Comments: Published by Geometry and Topology at}
\immediate\write\gtoutfile{\s\s http://www.maths.warwick.ac.uk/gt/GTVol\thevolumenumber/paper\thepapernumber.abs.html}
\immediate\write\gtoutfile{\noexpand\\}
\immediate\write\gtoutfile{}
\ifx\theasciiabstract\relax
\immediate\write\gtoutfile{\theabstract}\else
\immediate\write\gtoutfile{\theasciiabstract}\fi
\immediate\write\gtoutfile{}
\immediate\write\gtoutfile{\noexpand\\}
\immediate\write\gtoutfile{}
\immediate\closeout\gtoutfile}}  %%% end of definition of \makeheadfile

\def\maketitlepage{\maketitlep\makeheadfile}
\let\maketitle\maketitlepage

\def\ifplaintex{\expandafter\ifx\csname documentclass\endcsname\relax}

\ifplaintex 
\else
\fi

\expandafter\ifx\csname beginpicture\endcsname\relax
\expandafter\ifx\csname documentclass\endcsname\relax
\input pictex \else
\input prepictex \input pictex \input postpictex \fi\fi

\def\gt{{\mathsurround=0pt\it $\cal G\mskip-2mu$eometry \&\ 
$\cal T\!\!$opology}}        %  journal title in recommended style

\def\gtp{{\mathsurround=0pt\it $\cal G\mskip-2mu$eometry \&\ 
$\cal T\!\!$opology $\cal P\!$ublications}}  % GT publications

\def\lognumber#1{\def\thelognumber{#1}}
\def\volumenumber#1{\def\thevolumenumber{#1}}
\def\papernumber#1{\def\thepapernumber{#1}}
\def\volumeyear#1{\def\thevolumeyear{#1}}

\def\pagenumbers#1#2{\def\startpage{#1}\def\finishpage{#2}}
\def\published#1{\def\publishdate{#1}}
\def\proposed#1{\def\theproposer{#1}}
\def\seconded#1{\def\theseconders{#1}}
\def\received#1{\def\receiveddate{#1}}

\def\accepted#1{\def\accepteddate{#1}}
\def\asciititle#1{\def\theasciititle{#1}}

\long\def\asciiabstract#1{\long\def\theasciiabstract{#1}}
\def\asciikeywords#1{\def\theasciikeywords{#1}}

\let\\\par\let\thelognumber\relax
\let\thevolumenumber\relax\let\thepapernumber\relax
\let\thevolumeyear\relax\let\thesamplenumber\relax\let\startpage\relax
\let\finishpage\relax\let\publishdate\relax\let\receiveddate\relax
\let\reviseddate\relax\let\accepteddate\relax\let\theasciititle\relax
\let\theasciiauthors\relax
\let\theasciiabstract\relax\let\theasciikeywords\relax
\let\theasciiemail\relax\let\theshortauthors\relax\let\theshorttitle\relax

\long\def\maketitlep{   % start of definition of \maketitlep

\count0=\startpage

\gt\hfill      %   Journal title (top left) 
\beginpicture
\setcoordinatesystem units <0.33truein, 0.33truein> point at 2.2 0.9
\setplotsymbol ({$\cal G$})
\plotsymbolspacing=9truept
\circulararc 315 degrees from 0 1 center at 0 0
\setplotsymbol ({$\cal T$})
\circulararc 315 degrees from 1 -1 center at 1 0
\endpicture
\break
{\small\ifx\thesamplenumber\relax % sample?  
Volume \else Sample
\fi\thevolumenumber\ (\thevolumeyear)
\startpage--\finishpage\nl
Published: \publishdate}
\vglue 0.5truein plus 0.4fil minus 0.1truein

{\parskip=0pt\leftskip 0pt plus 1fil\def\\{\par\smallskip}{\ifplaintex\large
\else\Large\fi\bf\thetitle}\par\medskip}   

\vglue 0pt plus 0.1fil 

{\parskip=0pt\leftskip 0pt plus 1fil\def\\{\par}{\sc\theauthors}
\par\medskip}

\vglue 0pt plus 0.1fil 

{\small\parskip=0pt\let\newline\\
{\leftskip 0pt plus 1fil\def\\{\par}{\sl\theaddress}\par}
\expandafter\ifx\theemail\relax    % email address?
\relax\else\vglue 5pt plus 0.02fil minus 2pt\def\\{\stdspace{\rm 
and}\stdspace} 
\cl{Email:\stdspace\tt\theemail}\fi
\ifx\theurl\relax                  % URL given?
\relax\else\vglue 5pt plus 0.02fil minus 2pt\def\\{\stdspace{\rm 
and}\stdspace}
\cl{URL:\stdspace\tt\theurl}\fi\par}

\vglue 7pt plus 0.3fil minus 3pt

{\bf Abstract}
\vglue 5pt plus 0.1fil minus 2pt

\theabstract

\vglue 7pt plus 0.3fil minus 3pt

{\bf AMS Classification numbers}\quad Primary:\quad \theprimaryclass

Secondary:\quad \thesecondaryclass

\vglue 5pt plus 0.3fil minus 2pt

{\bf Keywords}\quad \thekeywords

\vglue 10pt plus 0.5fil minus 5pt

{\small  Proposed: \theproposer\hfill Received: \receiveddate\nl
Seconded: \theseconders\hfill 
\ifx\reviseddate\relax                         % paper revised?
Accepted: \accepteddate                        % no
\else
Revised: \reviseddate                          % yes
\fi}
\eject
}       %  end of definition of \maketitlep

\let\maketitlepage\maketitlep
\let\maketitle\maketitlepage

\font\phead=cmsl9 scaled 950
\font=cmsl9 scaled 1050
\font\pnum=cmbx10 scaled 913
\font\lnum=cmbx10 
\font\pfoot=cmsl9 scaled 950
\font=cmsl9 scaled 1050
\ifplaintex
\headline{\vbox to 0pt{\vskip -4.5mm\line{\small\phead\ifnum
\count0=\startpage ISSN 1364-0380 (on line)
1465-3060 (printed) \hfill {\pnum\folio}\else\ifodd\count0\def\\{ }% 
\ifx\theshorttitle\relax\thetitle\else\theshorttitle\fi\hfill{\pnum\folio}
\else\def\\{ and }{\pnum\folio}\hfill\ifx\theshortauthors\relax\theauthors
\else\theshortauthors\fi\fi\fi}\vss}}
\footline{\vbox to 0pt{\vglue 0mm\line{\small\pfoot\ifnum\count0=\startpage
\copyright\ \gtp\hfill\else
\gt, Volume \thevolumenumber\ (\thevolumeyear)\hfill\fi}\vss
}}
\else
\makeatletter
\def\@oddhead{{\small\lhead\ifnum\count0=\startpage ISSN 1364-0380 (on line)
1465-3060 (printed) \hfill {\lnum\number\count0}\else\ifodd\count0
\def\\{ }\ifx\theshorttitle\relax \thetitle \else\theshorttitle\fi\hfill
{\lnum\number\count0}\else\def\\{ and }{\lnum\number\count0}
\hfill\ifx\theshortauthors\relax 
\theauthors\else\theshortauthors\fi\fi\fi}}\def\@evenhead{\@oddhead}
\def\@oddfoot{\small\lfoot\ifnum\count0=\startpage\copyright\ \gtp\hfill\else
\gt, Volume \thevolumenumber\ (\thevolumeyear)\hfill\fi}
\def\@evenfoot{\@oddfoot}
\makeatother
\fi

\newwrite\gtoutfile
\long\gdef\makeheadfile{  %%% start of definition of \makeheadfile
{\def\\{, }\def\s{ }
\immediate\openout\gtoutfile head.xxx
\immediate\write\gtoutfile{To: math@arxiv.org}
\immediate\write\gtoutfile{Subject: put or rep NNNNN:pppp}
\immediate\write\gtoutfile{--text follows this line--}
\immediate\write\gtoutfile{Proxy-for: \ifx\theasciiauthors\relax
\theauthors\else\theasciiauthors\fi\s<\ifx\theasciiemail\relax\theemail\else\theasciiemail\fi>}
\immediate\write\gtoutfile{\noexpand\\}
\immediate\write\gtoutfile{Authors: \ifx\theasciiauthors\relax
\theauthors\else\theasciiauthors\fi}
{\def\\{ }\immediate\write\gtoutfile{Title: \ifx\theasciititle\relax
\thetitle\else\theasciititle\fi}}
\immediate\write\gtoutfile{Subj-class: GT or SG or MG etc}
\immediate\write\gtoutfile{MSC-class: \theprimaryclass\ifx\thesecondaryclass\relax\else, \thesecondaryclass\fi}
\immediate\write\gtoutfile{Journal-ref: Geom. Topol. \thevolumenumber
(\thevolumeyear) \startpage-\finishpage}
\immediate\write\gtoutfile{Comments: Published by Geometry and Topology at}
\immediate\write\gtoutfile{\s\s http://www.maths.warwick.ac.uk/gt/GTVol\thevolumenumber/paper\thepapernumber.abs.html}
\immediate\write\gtoutfile{\noexpand\\}
\immediate\write\gtoutfile{}
\ifx\theasciiabstract\relax
\immediate\write\gtoutfile{\theabstract}\else
\immediate\write\gtoutfile{\theasciiabstract}\fi
\immediate\write\gtoutfile{}
\immediate\write\gtoutfile{\noexpand\\}
\immediate\write\gtoutfile{}
\immediate\closeout\gtoutfile}}  %%% end of definition of \makeheadfile

\def\maketitlepage{\maketitlep\makeheadfile}
\let\maketitle\maketitlepage

%% file: 2002-2.bbl
\begin{thebibliography}
\bibitem{B}
\textbf{S\,B Bradlow}, \emph{Vortices in Holomorphic Line Bundles over Closed
  {K}{\"a}hler Manifolds}, Comm. Math. Phys. 135 (1990) 1--17

\bibitem{D1}
\textbf{S\,K Donaldson}, \emph{Monopoles, Knots, and Vortices}, (1997) lecture
  notes transcribed by Ivan Smith

\bibitem{D2}
\textbf{S\,K Donaldson}, \emph{Topological Field Theories and Formulae of
  {C}asson and {M}eng-{T}aubes}, from: ``Proceedings of the {K}irbyfest'',
  Geometry and Topology Monographs, 2 (1999)

\bibitem{FS}
\textbf{R Fintushel{\rm,} R Stern}, \emph{Knots, Links, and 4-Manifolds}, Invent.
  Math. 134 (1998) 363--400

\bibitem{HL2}
\textbf{M Hutchings{\rm,} Y-J Lee}, \emph{Circle-Valued {M}orse theory and
  {R}eidemeister Torsion}, preprint (1997)

\bibitem{HL1}
\textbf{M Hutchings{\rm,} Y-J Lee}, \emph{Circle-Valued {M}orse Theory,
  {R}eidemeister Torsion, and {S}eiberg-{W}itten Invariants of 3-Manifolds},
  Topology, 38 (1999) 861--888

\bibitem{IP}
\textbf{E Ionel{\rm,} T\,H Parker}, \emph{{G}romov Invariants and Symplectic Maps},
  Math. Ann. 314 (1999) 127--158

\bibitem{JT}
\textbf{A Jaffe{\rm,} C\,H Taubes}, \emph{Vortices and Monopoles: Structure of
  Static Gauge Theories}, Birkh{\"a}user, Boston, Mass. (1980)

\bibitem{MD}
\textbf{I\,G Macdonald}, \emph{Symmetric Products of an Algebraic Curve},
  Topology 1 (1962) 319--343

\bibitem{MT}
\textbf{G Meng{\rm,} C\,H Taubes}, \emph{\underline{SW} = {M}ilnor Torsion}, Math.
  Res. Lett. 3 (1996) 661--674

\bibitem{Milnor}
\textbf{J Milnor}, \emph{A Duality Theorem for {R}eidemeister Torsion}, Ann.
  Math. 76 (1962) 137--147

\bibitem{milnor2}
\textbf{J Milnor}, \emph{Lectures on the h-Cobordism Theorem}, Princeton
  Unversity Press (1965)

\bibitem{MOY}
\textbf{T Mrowka{\rm,} P Ozsv\'ath{\rm,} B Yu}, \emph{{S}eiberg-{W}itten Monopoles on
  {S}eifert Fibered Spaces}, Comm. Anal. Geom. 5 (1997) 685--791

\bibitem{OS1} \textbf{P Ozsv\'ath{\rm,} Z Szab\'o}, \emph{Holomorphic
Disks and Invariants for Rational Homology Three-Spheres}, 
{\tt arxiv:math.SG/0101206}

\bibitem{OS2} \textbf{P Ozsv\'ath{\rm,} Z Szab\'o}, \emph{Holomorphic
Disks and Three-Manifold Invariants: Properties and Applications},
{\tt arxiv:math.SG/0105202}

\bibitem{S}
\textbf{D Salamon}, \emph{{S}eiberg-{W}itten Invariants of Mapping Tori,
  Symplectic Fixed Points, and {L}efschetz Numbers}, from: ``Proceedings of 6th
  G\"okova Geometry-Topology Conference'' (1998)  117--143

\bibitem{Turaev1}
\textbf{V\,G Turaev}, \emph{{R}eidemeister Torsion in Knot Theory}, Russian
  Math. Surveys, 41 (1986) 119--182

\end{thebibliography}
